              \def\version{25 May 2020}	        	%

\documentclass[reqno,11pt]{amsart} 

\usepackage[T1]{fontenc}
\usepackage[utf8]{inputenc}
\usepackage{amsmath} 
\usepackage[mathscr]{eucal}
\usepackage{amssymb,bbm}
\usepackage{srcltx} 
\usepackage{dsfont}
\usepackage{hyperref}
\usepackage{color}
\usepackage{enumerate}
\usepackage{tikz}
\usepackage{boldline}
\usepackage{comment}
\usepackage{lscape}
\usepackage{graphicx}

\numberwithin{equation}{section}
 
\def\la{\lambda}

\def\lcvp{\lambda_{\rm vp}}
\def\lcve{\lambda_{\rm v}}
\def\lcdp{\lambda_{\rm dp}}
\def\lcde{\lambda_{\rm d}}
\def\lcnp{\lambda_{\rm np}}
\def\lcne{\lambda_{\rm n}}

\def\emptyset{\varnothing} 

\def\a{\alpha}

\def\d{{\rm d}} 
\def\e{\varepsilon} 
\def\l{l} 

\def\L{\Lambda}

\newfam\Bbbfam 
\font\tenBbb=msbm10 
\font\sevenBbb=msbm7 
\font\fiveBbb=msbm5 
\textfont\Bbbfam=\tenBbb 
\scriptfont\Bbbfam=\sevenBbb 
\scriptscriptfont\Bbbfam=\fiveBbb 
\def\Bbb{\fam\Bbbfam \tenBbb}

\def\one{\mathds{1}}

\newcommand{\C}     {\mathcal{C}} 
\newcommand{\R}     {\mathbb{R}} 
\newcommand{\Z}     {\mathbb{Z}} 
\newcommand{\N}     {\mathbb{N}} 
\newcommand{\M}     {\mathbb{M}} 
\renewcommand{\P}   {\mathbb{P}} 
 
\newcommand{\E}     {\mathbb{E}} 
\newcommand{\Q}     {\mathbb{Q}}

\def\1{{\mathchoice {1\mskip-4mu\mathrm l}      
{1\mskip-4mu\mathrm l} 
{1\mskip-4.5mu\mathrm l} {1\mskip-5mu\mathrm l}}} 
 
\def\comment#1{} 
\newtheoremstyle{thm}{2ex}{2ex}{\itshape\rmfamily}{} 
{\bfseries\rmfamily}{}{1.7ex}{} 
 
\newtheoremstyle{rem}{1.3ex}{1.3ex}{\rmfamily}{} 
{\itshape\rmfamily}{}{1.5ex}{}

\renewcommand{\theequation}{\thesection.\arabic{equation}} 
 
\newtheorem{theorem}{Theorem}[section] 
\newtheorem{lemma}[theorem]{Lemma} 
\newtheorem{prop}[theorem] {Proposition} 
\newtheorem{cor}[theorem]  {Corollary}

\theoremstyle{definition}
\newtheorem{defn}[theorem] {Definition} 
\newtheorem{example}[theorem] {Example}

\renewcommand{\d}{{\rm d}} 
 
\newcommand{\eps}{\varepsilon}

\newcommand{\esssup}{{\operatorname {esssup}}} 
\newcommand{\essinf}{{\operatorname {essinf}}}
\newcommand{\supp}{{\operatorname {supp}}} 
 
\newcommand{\dist}{{\operatorname {dist}}} 
\newcommand{\diam}{{\operatorname {diam}}}

\def\supp{\mathrm{supp}}

\newcommand\numberthis{\addtocounter{equation}{1}\tag{\theequation}}
\renewcommand{\e}   {{\operatorname e }}

\definecolor{Red}{rgb}{1,0,0}

 \title[Phase transitions for the Boolean model for Cox point processes]{Phase transitions for the Boolean model of continuum percolation for Cox point processes}
 \author[Elie Cali, Benedikt Jahnel, Andr\'as T\'obi\'as]{}
 
\begin{document}
 \maketitle
 \centerline{{\sc Benedikt Jahnel\footnote{Weierstrass Institute Berlin, Mohrenstra\ss e 39, 10117 Berlin, Germany, \texttt{Jahnel@wias-berlin.de}}, Andr\'as T\'obi\'as\footnote{TU Berlin, Stra\ss e des 17.~Juni 136, 10623 Berlin, Germany, \texttt{Tobias@math.tu-berlin.de}} and Elie Cali\footnote{Orange SA, 44 Avenue de la R\'epublique, 92326 Ch\^atillon, France, \texttt{Elie.Cali@orange.com}}}}
\renewcommand{\thefootnote}{}

\bigskip

\centerline{\small(\version)} 
\vspace{.5cm} 
 
\begin{quote} 
{\small {\bf Abstract:}} We consider the Boolean model with random radii based on Cox point processes. Under a condition of stabilization for the random environment, we establish existence and non-existence of subcritical regimes for the size of the cluster at the origin in terms of volume, diameter and number of points. Further, we prove uniqueness of the infinite cluster for sufficiently connected environments.  
\end{quote}




\setcounter{tocdepth}{3}

\section{Introduction and previous work}
The study of Boolean models based on stationary point processes in $\R^d$ traces back to the early 60's of the past century, when Gilbert~\cite{G61} introduced the model in the context of ad-hoc communication networks. The occurrence or absence of an infinite connected component in the set of points in $\R^d$ that have at least one point of the point process in their vicinity has served as a prototypical example of a phase transition ever since. More precisely, the model is defined via a {\em simple stationary point process} $X=(X_i)_{i \in I}$ in $\R^d$ with finite intensity $\lambda>0$ for $d \geq 1$. Every $X_i\in X$ carries an i.i.d.~mark $\varrho_i\ge 0$, which represents the {\em random interaction radius} of $X_i$. We will assume in the whole manuscript that $\P(\varrho>0)>0$, where $\varrho$ denotes a nonnegative radius random variable. The associated \emph{Boolean model} is then given by the union of balls centered at the points in $X$, with corresponding radii, i.e., 
\begin{align*}
\C=\bigcup_{i \in I} B_{\varrho_i}(X_i). 
\end{align*}
In the past 60 years this model has attracted much attention, especially when the underlying point process is given by a homogeneous Poisson point process with intensity $\lambda>0$, see for example~\cite{H85,MR96,G08,ATT16}, and many more contributions. In the following theorem, we start by summarizing some known results on coverage and uniqueness. 
\begin{theorem}[{\cite[Propositions 3.1 and 7.3, and Theorems 3.1, 3.6 and 7.4]{MR96}}]\label{Thm_1}
Let $\C$ be the Boolean model based on a stationary point process $X$ in $\R^d$ with finite intensity $\lambda>0$ and with i.i.d.~radii distributed as $\varrho$.   
\begin{enumerate}
    \item\label{co_co} \emph{Complete coverage.} Let $d \geq 1$. If $E[\varrho^d]=\infty$, then $\P(\C=\R^d)=1$. If $X$ is a homogeneous Poisson point process, then $\E[\varrho^d]<\infty$ implies that $\P(\C=\R^d)=0$.
        \item\label{d=1} \emph{One-dimensional triviality.} Let $d=1$ and $X$ a homogeneous Poisson point process. If  
     $\E[\varrho^d]<\infty$, then $\P(\C\text{ contains no unbounded connected component})=1$. 
     \item\label{Uniqu} \emph{Uniqueness.} Let $d\ge1$. If $\esssup(\varrho)=\infty$, then $\P(\C$ contains at most one unbounded component$)=1$. If $X$ is a homogeneous Poisson point process, then 
$\P(\C$ contains at most one unbounded component$)=1$ holds independent of the distribution of $\varrho$. 
\end{enumerate}
\end{theorem}

As already mentioned in Part~\eqref{d=1} of Theorem~\ref{Thm_1}, among the main objects of interest in the Boolean model are the maximal connected components, the {\em clusters}, and in particular, the {\em cluster that contains the origin}, which we denote by $C_o\subset\C$. 
For any measurable $A\subset \R^d$, writing $|A|$ for the Lebesgue volume, $\diam(A)=\sup\{|x-y|\colon x,y\in A\}$ for the diameter, and $X(A)=\#(X\cap A)$ for the number of points of $X$ in $A$, the following 
critical intensities can be defined via $C_o$,
\begin{center}
\begin{tabular}{lll} 
\hspace{-0.1cm}$\lcvp=\inf\{\lambda>0\colon \P(|C_o|=\infty)>0\}$& \text{and  for all }$s>0$& $\lcve(s)=\inf\{\lambda>0\colon \E[|C_o|^s]=\infty\}$,\\ 
\hspace{-0.1cm}$\lcdp=\inf\{\lambda>0\colon \P(\diam(C_o)=\infty)>0\}$& \text{and  for all }$s>0$& $\lcde(s)=\inf\{\lambda>0\colon \E[\diam(C_o)^s]=\infty\}$,  \\ 
\hspace{-0.1cm}$\lcnp=\inf\{\lambda>0\colon \P(X(C_o)=\infty)>0\}$& \text{and for all }$s>0$& $\lcne(s)=\inf\{\lambda>0\colon \E[X(C_o)^s]=\infty\}$. \\ 
\end{tabular}
\end{center}

In the next theorem, we present a summary of the known results on the non-triviality of these critical intensities in the case where $X$ is a homogeneous Poisson point process. 
\begin{theorem}[{\cite[Theorem 3.2]{MR96},~\cite[Theorems 2.1 and 2.2]{G08} and~\cite[Theorem 2]{GT18}}]\label{Thm_2} Let $d \geq 2$, $s>0$ and consider the Boolean model based on a homogeneous Poisson point process with i.i.d.~radii distributed as $\varrho$.  
\begin{enumerate}
\item\label{T2_SupCrit}We have that $\lcvp<\infty$, $\lcdp<\infty$, and $\lcnp<\infty$.
\item\label{T2_SubCrit} If $\E[\varrho^{d}] < \infty$, then $\lcvp>0$, $\lcdp>0$, and $\lcnp>0$.
\item\label{T2_SubMom} If $\E[\varrho^{d+s}] < \infty$, then $\lcve(s/d)>0$, $\lcde(s)>0$, and $\lcne(s/d)>0$.
\item\label{T2_NonSubMom} If $\E[\varrho^{d+s}] = \infty$, then $\lcve(s/d)=0$, $\lcde(s)=0$, and $\lcne(s/d)=0$.
\end{enumerate}
\end{theorem}
\begin{proof}
For Part~\eqref{T2_SupCrit}, $\lcvp<\infty$ is proved in~\cite[Theorem 3.2]{MR96} and $\lcdp<\infty$ is an immediate consequence of $\lcvp<\infty$. Further, $\lcnp<\infty$ is an immediate consequence of $\lcvp<\infty$ in case $\esssup(\varrho)<\infty$. However, via a straightforward coupling argument, it also holds in case $\esssup(\varrho)=\infty$.
For Part~\eqref{T2_SubCrit}, the fact that $\E[\varrho^{d}] < \infty$ implies $\lcvp>0$ is proved in~\cite[Theorem 2.1]{G08}. 
The statement that $\E[\varrho^{d}]<\infty$ implies $\lcdp>0$ follows from a straightforward adaptation of the proof of~\cite[Theorem 2.1]{G08}, noting that $2\sup_{x\in C_o}|x|\ge \diam (C_o)$.
Finally, the fact that $\E[\varrho^{d}]<\infty$ implies $\lcnp>0$ follows by yet another simple adaptation of the proof of~\cite[Theorem 2.1]{G08}, since 
\begin{align*}
\P(X(C_o)=\infty)\le \P(X(C_o)\ge X(B_\a(o))),
\end{align*}
where the right-hand side tends to zero for small $\lambda$ as $\a$ tends to infinity.
For parts~\eqref{T2_SubMom} and~\eqref{T2_NonSubMom}, note that the equivalence between $\lcve(s/d)>0$ and $\E[\varrho^{d+s}]$, as well as the equivalence between $\lcne(s/d)>0$ and $\E[\varrho^{d+s}]$, is proved in~\cite[Theorem 2]{GT18}. Finally, the equivalence between $\lcde(s)>0$ and $\E[\varrho^{d+s}]$ is provided in~\cite[Theorem 2.2]{G08}.
\end{proof}

In view of Theorem~\ref{Thm_2}, in particular we have that, as soon as $\E[\varrho^d]<\infty$, there exists a nontrivial  subcritical phase for volume percolation in the Poisson case, see Part~\eqref{T2_SubCrit}. Moreover, there are regimes where there exists a subcritical phase for volume percolation, however, the expected cluster size at the origin is infinite in terms of the number of points in the cluster. Similarly, various regimes for finite and infinite expected cluster size and expected diameter are available. In case of almost surely bounded radii, we have a number of identities for the critical intensities in the case of Poisson point processes, see~\cite[Theorem 3.4, 3.5, and  Proposition 3.2]{MR96}, a topic that we will not address further in this manuscript. 


\bigskip
The vast majority of available results for continuum percolation are for the Boolean model based on Poisson point processes. However, apart from the statements for general stationary point processes already presented above, there are by now a number of assertions available for continuum percolation based on a variety of point processes other than the Poisson point process. For example, in~\cite{M75,S13,J16,M18} percolation properties for the Boolean model based on a large class of stationary Gibbs point processes are considered for fixed radii. In~\cite{CD14}, percolation in the Boolean model is proved to exist for dense Gibbs point processes with quermass interaction, where the radii are random but not i.i.d.. Further, in~\cite{GKP16}, the authors study a class of repelling point processes in $\R^2$ that include the Ginibre ensemble as well as the Gaussian zero process, with fixed radii. For the Gaussian zero process, they establish uniqueness of the infinite cluster and non-triviality of continuum percolation. For the Ginibre ensemble, they also prove uniqueness of the infinite cluster; the non-triviality of continuum percolation for this point process was already verified in \cite{BY14}.
In the same direction, in~\cite{BY13}, the authors study clustering and percolation for a family of stationary point processes that exhibit negative association, such as general determinantal point processes and some perturbed lattices, see also~\cite{GY05}. These processes may cluster less than the Poisson point process but still have non-trivial percolation properties in the associated Boolean model with non-random radii. The paper~\cite{BY13} also presents an example of a stationary Cox point process with stronger clustering properties than the Poisson point process, but which nevertheless percolates for all positive non-random radii. 
In~\cite{HJC19}, non-triviality of percolation is shown for the Boolean model with non-random radii, based on stabilizing Cox point processes with sufficiently connected support. This was the starting point for our investigation presented in this paper, and we will refer to these results frequently in what follows. Finally, in~\cite{G09}, conditions are presented for the existence of subcritical regimes for percolation for Boolean models based on stationary point processes. The derivation of these conditions is very similar to our approach, and in fact, checking these conditions provides an alternative strategy for the proof of some of our results.

\section{Setting and main results}
\begin{figure}[!htpb]
\centering
\input{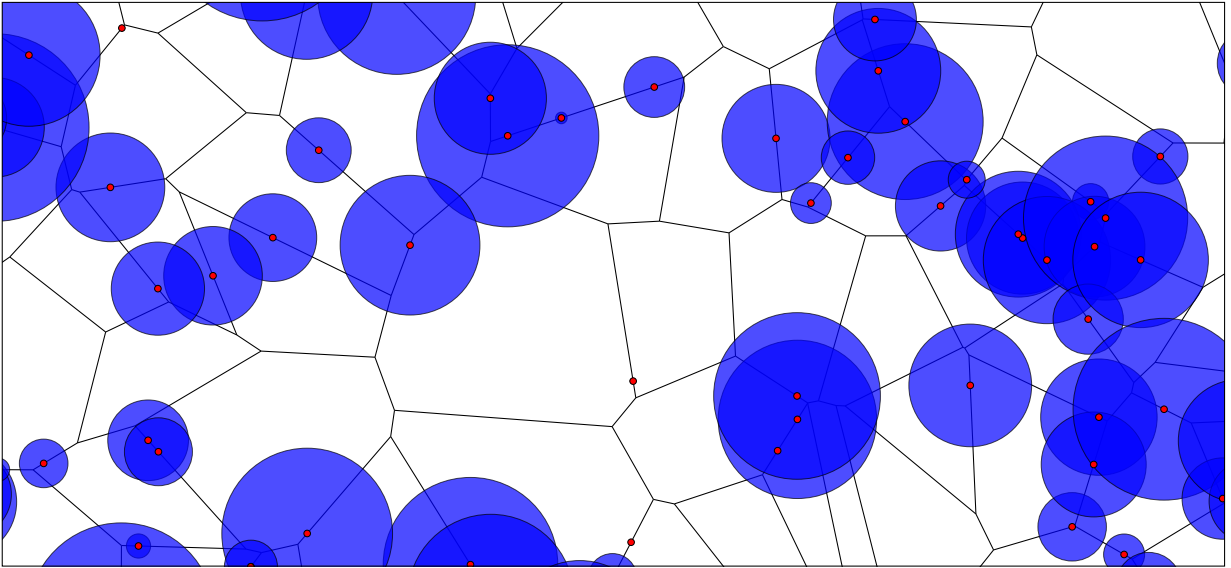}
\caption{Realization of a Boolean model (blue) based on a Cox point processes (red) with directing measure given by a realization of a Poisson--Voronoi tessellation (black).}
\label{Fig-CoxPVT}
\end{figure}

We are interested in the setting where $X$ is a Cox point process, i.e., a Poisson point process with a random intensity measure. To be more precise, we consider random elements $\L$ in the space $\M$ of Borel measures on $\R^d$ equipped with the evaluation $\sigma$-algebra. In the whole manuscript we assume $\L$ to be stationary and normalized so that it satisfies $\E[\Lambda([0,1]^d)]=1$. Then, for $\lambda>0$, we let $X$ be a {\em Cox point process} with intensity measure $\lambda \Lambda$ and call $\L$ its {\em directing measure}. The particular case where $\Lambda$ equals the Lebesgue measure corresponds to a homogeneous Poisson point process with intensity $\lambda$.  Figure~\ref{Fig-CoxPVT} provides an illustration of a realization of the Boolean model based on a Cox point process, where the directing measure $\Lambda$ is given by the edge-length measure on a Poisson--Voronoi tessellation, see Section~\ref{sec-SingularExamples} for more details. 

Under conditions of stabilization and connectedness, see below, we essentially reproduce the known picture for the Poisson case also for the Cox case. Let us start by noting that some known results for general stationary point processes, as for example the one-dimensional triviality or the first part of the complete-coverage result from Theorem~\ref{Thm_1}, immediately transfer to the case of Cox point processes. In the next few results, we recover the statements of Theorem~\ref{Thm_1} related to Poisson point processes. We start by giving an extension of the second statement in Part~\eqref{co_co} in Theorem~\ref{Thm_1}. Here, $\L$ is called {\em ergodic} if $\P(\L\in E)\in\{0,1\}$ for all measurable and translation-invariant events $E\subset \M$. 
\begin{prop}[No complete coverage]\label{Prop_1}
Let $d \geq 1$ and $\C$ be the Boolean model based on a Cox point process with stationary intensity measure $\lambda\Lambda$ and with i.i.d.~radii distributed as $\varrho$. If $\E[\varrho^d]<\infty$, then $\P(\C=\R^d)<1$ for all $\lambda>0$. If additionally $\Lambda$ is ergodic, then $\P(\C=\R^d)=0$. 
\end{prop}
We next lift the result in Part~\eqref{d=1} in Theorem~\ref{Thm_1} to the level of Cox point processes with ergodic intensity measures. 
\begin{prop}[One-dimensional triviality]\label{Prop_2}
Let $d=1$ and $\C$ be the Boolean model based on a Cox point process with ergodic intensity measure $\lambda\Lambda$ and with i.i.d.~radii distributed as $\varrho$. If  $\E[\varrho^d]<\infty$, then $\C$ contains no unbounded component, almost surely.
\end{prop}
Let us note that (at least) for $\Lambda$ a convex combination of finitely many ergodic measures, the statement of Proposition~\ref{Prop_2} should also hold. For the uniqueness result, we need the following condition on the connectedness of the directing measure $\Lambda$. 
Let us write $Q_r(x)=[-r,r]^d+x$ for the box with sidelength $2r$, centered at $x\in\R^d$ and abbreviate $Q_r=Q_r(o)$ and $A^{\rm c}(x)=\R^d\setminus A$ for $A \subseteq \R^d$. Further, let us write 
$$\supp(\L)=\{x\in \R^d\colon \L(Q_\eps(x))>0\text{ for every }\eps>0\}$$
for the support of a measure $\L$ and $\L_A$ for the restriction of $\L$ to the set $A$, i.e., $\L_A(B)=\L(A\cap B)$ for measurable $A,B\subset\R^d$. 
\begin{defn}[Essential $r$-connectedness]\label{Conn}
Let $r>0$. We say that the stationary random measure $\L$ is \emph{essentially $r$-connected} if there exists a random field of \textit{connectivity radii} $R=\{R_x\}_{x\in\R^d}$, defined on the same probability space as $\L$, such that
\begin{enumerate}
\item $(\L, R)$ are jointly stationary, 
\item $\lim_{\a\uparrow\infty}\P(\sup_{y\in Q_{\a} \cap\,  \Q^d}R_y \ge \a)=0$, and
\item for all $\alpha \geq 1$, whenever $\sup_{y\in Q_{2\a} \cap\,  \Q^d}R_y < \alpha/2$, we have that for all $x,y\in \supp(\L_{Q_\a})$ there exists a finite sequence of points $(x_1,\dots, x_l)\subset\supp(\L_{Q_{2\a}})$ such that $|x_i-x_{i+1}|<r$ for all $i\in\{0,1,\dots, l+1\}$ where $x=x_0$ and $y=x_{l+1}$.   
\end{enumerate} 
\end{defn}
Let us mention that our definition of essential $r$-connectedness is different from the definition of {\em asymptotic essential connectedness} as presented in~\cite{HJC19}. One of the differences is that essential $r$-connectedness does not imply stabilization, see below. However, for many directing measures both connectedness conditions can be verified. We comment on this further and present examples for essentially $r$-connected directing measures in Section~\ref{Sec-Examples}. 
Equipped with these definitions, we can now state our result on the uniqueness of the unbounded component for Cox point processes. 
\begin{theorem}[Uniqueness]\label{Prop_3}
Let $d \geq 1$ and $\C$ be the Boolean model based on an ergodic Cox point process with stationary intensity measure $\lambda\Lambda$ and with i.i.d.~radii distributed as $\varrho$. If $r=\esssup(\varrho)<\infty$ and $\Lambda$ is essentially $r$-connected, then $\C$ has at most one unbounded cluster, almost surely. 
\end{theorem}
Recall that if $\esssup(\varrho)=\infty$, then Theorem~\ref{Thm_1} already implies uniqueness of the infinite cluster for any stationary point process. 
Note that the Propositions~\ref{Prop_1},~\ref{Prop_2} and Theorem~\ref{Prop_3} establish a version of Theorem~\ref{Thm_1} for Cox point processes. Their proof can be found in Section~\ref{sec-Basicproofs}.
Before we come to our main result, in which we present conditions under which a version of Theorem~\ref{Thm_2} is available for Cox point processes, we need to introduce the following condition on spatial decorrelation of the directing measures. For this we define for all measurable $A,B\subset\R^d$ the distance $\dist(A,B)=\inf\{|x-y|\colon x\in A\text{ and }y\in B\}$. 
\begin{defn}[Stabilization]\label{Stab}
We say that the stationary random measure $\L$ is \emph{$\phi$-stabilizing} if there exists a random field of \textit{stabilization radii} $R=\{R_x\}_{x\in\R^d}$, defined on the same probability space as $\L$, such that
	\begin{enumerate}
    \item $(\L, R)$ are jointly stationary,
    \item $\lim_{\a\uparrow\infty}\phi(\a)=0$, where $\phi(\alpha)=\P(\sup_{y\in Q_\a \cap\,  \Q^d}R_y \ge \a)$, and
    \item for all $\a \ge 1$, the random variables
    $$\big(f(\L_{Q_\a(x)})\one\{\sup_{y \in Q_\a(x) \cap \, \Q^d}R_y < \a\}\big)_{x\in\psi}$$
        are independent for all bounded measurable functions $f\colon \M\to[0,\infty)$ and finite $\psi\subset\R^d$, as long as $\dist(x,\psi\setminus x)> 3\a$ for all $x\in\psi$. 
\end{enumerate}
\end{defn}
Note that $\phi$-stabilization for $\L$ implies that $\L$ is ergodic. This can be easily seen by verifying a mixing condition for local events.  We present examples for $\phi$-stabilizing directing measures in Section~\ref{Sec-Examples}. 
\begin{theorem}\label{Thm_3} Let $d \geq 2$, $s>0$ and consider the Boolean model based on a Cox point process with stationary intensity measure $\lambda\Lambda$ and with i.i.d.~radii distributed as $\varrho$.  
\begin{enumerate}
\item\label{T3_SupCrit} If $\L$ is
$\phi$-stabilizing with sufficiently large $\esssup(\varrho)$, then $\lcvp<\infty$, $\lcdp<\infty$ and $\lcnp<\infty$.
\item\label{T3_SubCrit} Let $\E[\varrho^{d}] < \infty$. If $\L$ is $\phi$-stabilizing, then $\lcvp>0$, $\lcdp>0$ and $\lcnp>0$.
\item\label{T3_SubMom} Let $\E[\varrho^{d+s}] < \infty$. If $\L$ is $\phi$-stabilizing with $\int_{0}^\infty\a^{s-1}\phi(\a)\d \a<\infty$, then 
$\lcve(s/d)>0$ and $\lcde(s)>0$.
Further, $\lcne(s/d)>0$ if $\L$ is $\phi$-stabilizing with $\int_{0}^\infty\a^{s-1}\phi(\a)\d \a<\infty$ and
        \begin{align}\label{Erg_cond}
    \int_{0}^\infty \a^{s-1}\P(\L(B_{\a})\ge c\a^d)\d\a<\infty\quad\text{for some }c>0. 
    \end{align}
\item\label{T3_NonSubMom} Let $\E[\varrho^{d+s}] = \infty$. If $\L$ is ergodic, then $\lcve(s/d)=0$, $\lcde(s)=0$ and $\lcne(s/d)=0$.
\end{enumerate}
\end{theorem}
Let us comment on the result, the proof of which can be found in Section~\ref{sec-Mainproofs}. A version of the statement that $\lcvp<\infty$ has been first proved under a condition of asymptotic essential connectedness in~\cite{HJC19} for fixed positive radii. 
The statement that $\lcvp<\infty$ for stabilizing directing measures if the (fixed) radius is sufficiently large originates from~\cite[Corollary 2.5]{T18}. The same assertion for random radii follows from \cite[Theorem 1.1]{JT19b}, but since the papers \cite{T18, JT19b} use a rather different notation,
we provide a self-contained proof for part \eqref{T3_SupCrit} of Theorem~\ref{Thm_3} in the present paper. Regarding Part~\eqref{T3_SubCrit} of Theorem~\ref{Thm_3}, the statement that $\lcvp>0$, for stabilizing directing measures, has been proven for the case of fixed nonnegative radii in~\cite{HJC19}. Moreover, the statements about $\lcvp>0$ and $\lcdp>0$ can be alternatively established by proving that certain general conditions provided in~\cite{G09} are applicable for stabilizing Cox point processes, see also our comments in Section~\ref{sec-Mainproofs}. 
This also holds for the statements about $\lcve(s/d)>0$ and $\lcde(s)>0$ in Part~\eqref{T3_SubMom} of Theorem~\ref{Thm_3}. In particular,~\cite[Theorem 1.3]{G09} provides the statements for $\lcvp>0$, $\lcdp>0$, $\lcve(s/d)>0$ and $\lcde(s)>0$ with respect to general stationary point processes that feature finite dependencies, in the spirit of our definition of $b$-dependence below (cf.~Section~\ref{ex-abs_enviroment}). Although checking the conditions of \cite{G09} is relatively straightforward, it is not very insightful. Hence, we decided to present a self-contained proof of these parts of Theorem~\ref{Thm_3} instead, extending the methods of \cite{G08} to the case of Cox point processes. The integrability conditions on the stabilization probability $\phi$ of the directing measures are often easy to verify for particular examples, as we will exhibit in Section~\ref{Sec-Examples}. Finally, Condition~\eqref{Erg_cond} establishes integrability of the deviation of the directing measure from its expectation. It can be guaranteed with the help of some moment conditions on the environment, which we present now and further use in Section~\ref{Sec-Examples}. 
\begin{lemma}\label{lem_4}
Condition~\eqref{Erg_cond} holds if any of the following two conditions holds,
\begin{align}
\limsup_{\a\uparrow\infty}\a^{-d}\log\E[\exp(\beta\L(B_\a))]&<\infty\qquad\text{for some }\beta>0,\label{eq:cond1}\\
\limsup_{\a\uparrow\infty}\a^{s-d\beta+\eps}\E\big[\big|\L(B_\a)-|B_\a|\big|^\beta\big]&<\infty\qquad\text{for some }\beta\ge 1\text{ and }\eps>0. \label{eq:cond2}
\end{align}
\end{lemma}

We prove this lemma in Section~\ref{sec-Exampleproofs}. In the next section we present and discuss a number of examples of Cox point processes in order to illustrate the generality and also limitations of our results.  
\section{Examples}\label{Sec-Examples}
\subsection{Boundedness, $\phi$-stabilization, essential $r$-connectedness and ergodicity}
As has been discussed already in~\cite{HJC19}, often, directing measures $\L$ can be categorized as either being absolutely continuous or singular with respect to the Lebesgue measure on $\R^d$. 

\subsubsection{Absolutely continuous directing measures}\label{ex-abs_enviroment}
We say that $\L$ is an {\em absolutely continuous directing measure} if $\L(\d x)=\ell_x\d x$ is given via a non-negative random field $(\ell_x)_{x\in \R^d}$, where $\d x$ denotes the Lebesgue measure on $\R^d$. A stationary and  absolutely continuous directing measure is called {\em bounded} if there exists $M>0$ such that $\ell_o<M$ almost surely. Let us consider some examples of such random fields.

\begin{example}[Bounded absolutely continuous directing measures]\label{Ex_1}
Bounded environments can be constructed via random fields of the form $\ell_x=\lambda_1\one\{x\in \Xi\}+\lambda_2\one\{x\not\in\Xi\}$, where $\Xi\subset\R^d$ is a random closed set and $\lambda_1,\lambda_2\ge 0$. For example, $\Xi$ could be given as an independent Boolean model $\bar \C$,
see Figure~\ref{Fig-CoxBool} for an illustration. 
\begin{figure}[!htpb]
\centering
\input{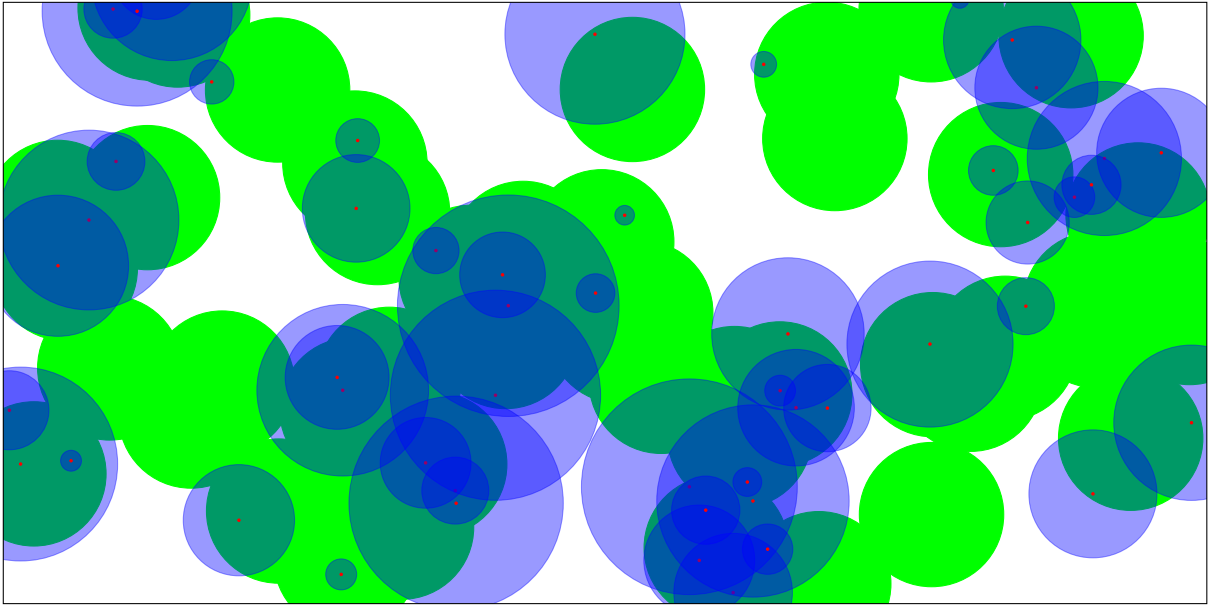}
\caption{Realization of a Boolean model (blue) based on a Cox point process (red) with directing measure given by an independent Poisson--Boolean model with fixed radii (green).}
\label{Fig-CoxBool}
\end{figure}
Note that in case both $\lambda_1,\lambda_2> 0$, then $\L$ is essentially $r$-connected for any $r>0$. In this case it is also asymptotically essentially connected in the sense of~\cite{HJC19}. 
\end{example}

As a corollary of Theorems~\ref{Thm_2} and \ref{Thm_3}, let us note that the existence of a subcritical phase for a Cox--Boolean model where the Cox point process has an absolutely continuous, bounded, and ergodic directing measure, behaves the same way as in the case of the Poisson point process. 
\begin{cor}\label{cor-bdd}
The statements of Part~\eqref{T3_SubCrit} of Theorem~\ref{Thm_2} remain true if Poisson point processes are replaced by Cox point processes with ergodic, absolutely continuous, and bounded directing measures.
\end{cor}
The proof of this corollary can be found in Section~\ref{sec-Exampleproofs}. The analogous statement does in general not hold for the existence of a supercritical phase as presented in Part~\eqref{T3_SupCrit} of Theorem~\ref{Thm_2}, as already mentioned in~\cite{HJC19}. For example, it may fail in the case of Example~\ref{Ex_1} when $\bar\C$ is the Poisson--Boolean model with constant radii and either $\la_1$ or $\la_2$ is zero, see~\cite{T18}.

We say that a $\phi$-stabilizing directing measure $\L$ is {\em $b$-dependent} if there exists $b>0$, such that $\phi(\a)=0$ for all $\a>b$. Further, we call a $\phi$-stabilizing directing measure $\L$ {\em exponentially stabilizing} if there exist $c>0$, such that $\phi(\a)\le\exp(-c\a)$ for all sufficiently large $\a>0$. 
\begin{example}[Unbounded $b$-dependent absolutely continuous directing measures]\label{Ex_2}
As an example of an unbounded absolutely continuous environment consider the {\em shot-noise field}, where $\ell_x=\sum_{i\in I}\kappa(Y_i-x)$, with $(Y_i)_{i\in I}$ a homogeneous Poisson point process and $\kappa\colon \R^d\to[0,\infty)$ integrable. If $\kappa$ is additionally compactly supported and the diameter of the support is given by $b$, then the associated directing measure is $b$-dependent. If the support of $\kappa$ has Lebesgue measure equal to infinity, then there exists no $\phi$ such that the associated directing measure is $\phi$-stabilizing. For any $r>0$, in general, the shot-noise field is also not essentially $r$-connected and also not asymptotically essentially connected in the sense of~\cite{HJC19}. However, there are cases when it is asymptotically essentially connected, see~\cite[Section 2.5.1]{T18}. 
\end{example}

\begin{example}[Unbounded stabilizing absolutely continuous directing measures]\label{Ex_3}
Keeping the Example~\ref{Ex_1} in mind, consider the directing measure based on the random field given by 
$\ell_x=\la\sum_{i\in I}\one\{|Y_i-x|<\varrho_i\}$
where $\la>0$ and $(Y_i)_{i\in I}$ is an independent homogeneous Poisson point processes with intensity $\mu$ and with i.i.d.~marks $\varrho_i$, distributed according to $\varrho$. This is essentially an environment coming from a Poisson--Boolean model with random radii, where instead of considering $\C$, we count the number of balls in $\C$ containing a given point $x$. In particular, this environment is unbounded. 
Again, if $b=\esssup(\varrho)<\infty$, then $\L$ is $b$-dependent and stochastically dominated from above by a shot-noise field, see Example~\ref{Ex_2}. If even $\essinf(\varrho)>0$, then it is also stochastically dominated from below by a shot-noise field.  For $\esssup(\varrho)<\infty$, this model is not essentially $r$-connected for any $r>0$. However, it is asymptotically essentially connected, in the sense of~\cite{HJC19}, for all fixed finite $\varrho$ and sufficiently large $\mu$. If $\esssup(\varrho)=\infty$ and $\E[\varrho^{d+s}] < \infty$ for $s>0$, then $\L$ is $\phi$-stabilizing with $\int_{0}^\infty\a^{s-1}\phi(\a)\d \a<\infty$. Indeed, define the stabilization radii 
$$R_y=\sup\{\varrho_i\colon |Y_i-y|<\varrho_i, Y_i\in Y\},$$
as the largest radius of a ball in $\C$ that contains $y\in\R^d$. Then, indeed, for $x,z\in\R^d$ with $|x-z|>3\a$, 
 \begin{align}\label{eq: stab}
 \L_{Q_\a(x)}\one\{\sup_{y \in Q_\a(x) \cap \, \Q^d}R_y < \a\}\qquad\text{ and }\qquad\L_{Q_\a(z)}\one\{\sup_{y \in Q_\a(z) \cap \, \Q^d}R_y < \a\}
 \end{align}
are independent, since, under the stabilization event~\eqref{eq: stab}, $\big(\C\cap Q_\a(x)\big)\cap B_{\varrho_i}(Y_i)=\emptyset$ for all $Y_i\in Y\cap Q^{\rm c}_{2\a}(x)$. Moreover, using Campbell's formula, 
\begin{align*}
\phi(\a)&\le \P\big(\#(Y_i\in Y\colon \varrho_i\ge \a,\, |Y_i|<\varrho_i+2\a)\ge 1\big)\le \mu\int_\a^\infty|B_{2\a+r}|\nu(\d r), 
\end{align*}
where $\nu=\P \circ \varrho^{-1}$, and hence for all $s>0$, there exists a constant $c$, only depending on $s$ and the dimension, such that 
\begin{align*}
\int_{0}^\infty\a^{s-1}\phi(\a)\d \a\le c\int_{0}^\infty r^{d+s}\nu(\d r),
\end{align*}
which is finite by assumption. 
Note that such an example, where there is stabilization but it is not necessarily exponential, was not provided in~\cite{HJC19,T18}. 
%
%
Finally, if for some $\beta>0$ we have that $\E[\exp(\beta\varrho)] < \infty$, then 
\begin{align*}
\phi(\a)&\le \mu\int_\a^\infty|B_{2\a+r}|\nu(\d r)\le \mu v_d\exp(-\beta\a)\int_0^\infty r^d\exp(\beta r/3)\nu(\d r)\le c \exp(-\beta\a),
\end{align*}
for some constant $c>0$, writing $v_d=|B_1|$, since $r^d\exp(\beta r/3)\le \exp(\beta r)$ for sufficiently large $r$. Hence, in this case, the directing measure is even exponentially stabilizing.  Similarly, one can show that if $\E[\varrho^{d}] <\infty$ but $\E[\varrho^{d+s}] =\infty$, then the associated $\L$ is $\phi$-stabilizing but $\int_{0}^\infty\a^{s-1}\phi(\a)\d \a=\infty$. Finally, note that in case $\E[\varrho^{d}] =\infty$, then $\ell_o$ is equal to infinity, almost surely, and hence the model is not well-defined, see~\cite[Section 3.1]{MR96}. 
\end{example}

\subsubsection{Singular directing measures}\label{sec-SingularExamples}
We say that a directing measure $\L$ is {\em singular} if $\L(\d x)$ is almost surely singular with respect to the Lebesgue measure $\d x$ on $\R^d$. We already gave an illustration of a singular directing measure given by the edge-length measure of a {\em Poisson--Voronoi tessellation} in Figure~\ref{Fig-CoxPVT}. More formally, a large class of singular environments can be constructed via {\em random segment processes} $S\subset\R^d$, by defining directing measures $\L(A)=\nu_1(S\cap A)$ for any measurable $A\subset\R^d$, where $\nu_1$ denotes the one-dimensional Hausdorff measure. 
\begin{example}[Stabilizing singular directing measures]\label{Ex_4}
A particularly interesting class of segment processes are tessellation processes in $\R^2$. Most prominently the Poisson--Voronoi and the {\em Poisson--Delaunay tessellations}, see~\cite{M89} for more background on tessellation processes. It has already been observed in~\cite{HJC19} that both the Poisson--Voronoi and the Poisson--Delaunay tessellation give rise to exponentially stabilizing directing measures that are asymptotically essentially connected, see~\cite[Section 3]{HJC19} for details. From this it is clear that they are also essentially $r$-connected for any $r>0$.  
\end{example}

\begin{example}[Essentially $r$-connected but not stabilizing directing measures]\label{Ex_5}
Another particularly interesting class of segment processes is the \emph{Manhattan grid} in $\R^2$, see for example~\cite{JT19}, which is defined as follows. Let $Y=(Y_{\rm v},Y_{\rm h})$ be the tuple where $Y_{\rm v}=\{Y_{i,\rm v}\}_{i\in I_{\rm v}}$ and $Y_{\rm h}=\{Y_{i,\rm h}\}_{i\in I_{\rm h}}$ are two independent simple stationary point processes on $\R$. 
Then, the Manhattan grid is defined as
\begin{align*}
S=S(Y)=\bigcup_{i\in I_{\rm v},\, j\in I_{\rm h}}(\R\times \{Y_{i,\rm h}\})\cup(\{Y_{j,\rm v}\}\times\R).
\end{align*}
The singular directing measure $\Lambda$ corresponding to the Manhattan grid is essentially $r$-connected for all $r>0$. Indeed, for $(x,y) \in \R^2$ define the connectivity radius
\[ R_{(x,y)} = \inf \big\{ r>0 \colon [y,y+r) \cap  Y_{\rm h} \neq \emptyset \text{ and } [x,x+r) \cap Y_{\rm v} \neq \emptyset \big\}. \]  
In words, $R_{(x,y)}$ is the smallest distance within which, measured from $(x,y)$, there is a horizontal line of the Manhattan grid above $(x,y)$ and a vertical line on the right of $(x,y)$. Clearly, for $R=(R_{(x,y)})_{(x,y) \in \R^2}$, points (1) and (2) of Definition~\ref{Conn} are satisfied. Assume now that for some $\alpha>0$, $\sup_{(x,y) \in Q_{\alpha} \cap \Q^d} R_{(x,y)}<\alpha$. Then in particular, for $(x,y)=(0,\alpha)$, there is a horizontal line $l_{\rm h}$ of the Manhattan grid above $(x,y)$ within distance at most $\alpha$. Further, for $(x,y)=(\alpha,0)$, there is a vertical line $l_{\rm v}$ on the right of $(x,y)$ within distance at most $\alpha$. Now, any vertical (respectively horizontal) line of $S$ intersecting with $Q_{\alpha}$ intersects with $l_{\rm h}$ (respectively $l_{\rm v}$) within $Q_{2n}$. Since $l_{\rm h}$ and $l_{\rm v}$ also intersect within $Q_{2n}$, it follows that $\supp(\Lambda_{Q_\alpha})$ is entirely connected within one single connected component of $\supp(\Lambda_{Q_{2\alpha}})$, which implies that condition (3) of Definition~\ref{Conn} also holds.

Since $S$ consists of infinite lines, $\Lambda$ is not $\phi$-stabilizing for any $\phi$. Hence, the vast majority of the results of \cite{HJC19} is not applicable in the case of the Cox point process with directing measure $\lambda\Lambda$. In particular, for constant radii, also in the simplest nontrivial case when both $Y_{\rm v}$ and $Y_{\rm h}$ are Poisson processes, both existence of a subcritical case and the one of a supercritical phase for percolation are unknown. Nevertheless, thanks to the fact that $\Lambda$ is essentially $r$-connected and also ergodic, some of the results of the present paper that are not necessarily true for any stationary Cox point process can be applied for Manhattan grids. Namely, Theorem~\ref{Prop_3} holds for them, i.e., there is at most one unbounded cluster, also if the radii are bounded. Further, Part~\eqref{T3_NonSubMom} of Theorem~\ref{Thm_3} implies that for $\varrho$ sufficiently heavy-tailed, if a certain moment of the volume, diameter and point-number of $C_o$ is infinite for the Poisson point process, then the same holds for the Manhattan grid. 

An even more fundamental segment process in $\R^2$ is the {\em Poisson line tessellation}, see for example~\cite{M89,JT19}, which is also ergodic but not $\phi$-stabilizing for any $\phi$. We expect that it is also essentially $r$-connected and this can be verified similarly to the case of Manhattan grids, but we refrain from presenting here the details. 
\end{example}

\subsubsection{Pathological behaviour of non-ergodic examples}

In the case of constant radii where $\varrho \equiv r>0$ holds almost surely, it is easy to show that $\lcvp=\lcdp=\lcnp$. Now, as already mentioned, if $\Lambda$ is $\phi$-stabilizing, then this critical intensity is positive for any $r>0$ and finite for all sufficiently large $r>0$. 
In Section~\ref{ex-abs_enviroment} we explained that $\lcvp=\lcdp=\lcnp=\infty$ occurs for some relevant $\phi$-stabilizing examples for small $r>0$. In contrast, the lack of a subcritical phase is often thought of as a pathology. A general construction of examples of $\Lambda$ with this property can be found in \cite[Section 4.2.3.3]{T19}. The simplest one among these examples is the mixed Poisson point process, which also lacks a subcritical phase in the case of general nonzero random radii, as the following example shows.

\begin{example}[Zero critical intensity]\label{ex-mixedPPP}
The stationary \emph{mixed Poisson point process} is the Cox point process with directing measure $\Lambda(\d x)=Z \d x$, where $Z$ is a non-negative random variable with $\E[Z]=1$. Then, $\Lambda$ is not ergodic and in particular not $\phi$-stabilizing. We claim that for $d \geq 2$, if $\esssup(Z)=\infty$ and $\varrho$ satisfies $\P(\varrho>0)>0$, then $\lcvp=\lcdp=\lcnp=0$. Indeed, then we can choose $\eps>0$ and $p \in (0,1]$ such that $\P(\varrho>\eps)=p>0$. Now, let $\lambda>0$ and consider the Cox point process $X$ with directing measure $\lambda Z\d x$ for $\lambda>0$. Then, the independent thinning $X^\eps=\{ X_i \in X \colon \varrho_i>\eps \}$ of $X$ is a Cox point process with nonzero intensity $Zp \lambda\d x$. Now, in the event $\{ Z>2\lambda_{\rm c}(\eps)/(p\lambda) \}$ that has positive probability, this intensity is greater than $2\lambda_{\rm c}(\eps)$, where $\lambda_{\rm c}(\eps)$ denotes the critical intensity for percolation of the Boolean model with constant radii $\eps$ based on a Poisson point process.  Thus, $\bigcup_{X_i \in X^\eps} B_{\varrho_i}(X_i)$ stochastically dominates the Boolean model of a Poisson point process of supercritical intensity, hence, almost surely it contains a cluster including infinitely many points of $X$, which also has infinite volume as well as infinite diameter. It follows that $\lcvp=\lcdp=\lcnp=0$, as wanted. Note that in case $\essinf(Z)>0$, the mixed Poisson point process is also essentially $r$-connected for any $r>0$. However, since $\Lambda$ is not ergodic, our uniqueness result is technically not applicable, but thanks to a comparison to the Poisson point process, there cannot be more than one unbounded cluster. Still,  the number of unbounded clusters may not be constant almost surely. 
\end{example}

In the following, similar example, $\lcvp=\lcdp=\lcnp=0$ also holds if $\esssup(\varrho)$ is sufficiently large but $\lcvp=\lcdp=\lcnp=\infty$ if $\esssup(\varrho)$ is small and positive.

\begin{example}[Zero critical intensity only for large radii]\label{ex-mixedCPP}
Let $\Lambda'$ be any $\phi$-stabilizing intensity measure such that $\lcnp(r)=\infty$ for sufficiently small constant radii $r>0$ and $\lcnp(r)<\infty$ for large enough constant radii $r$. Recall that for non-random radii $\lcnp(r)=\lcdp(r)=\lcvp(r)$. Let $r_{\rm c}=\inf\{r>0\colon \lcnp(r)<\infty\}$ denote the critical radius for the existence of a supercritical phase.
Let now $Z$ be an unbounded random variable with $\E[Z]=1$ that is independent of $\Lambda'$ and define the intensity measure $\Lambda$ via $\Lambda(\d x)=Z\Lambda'(\d x)$. Then, $\Lambda$ is not ergodic. Now, if $r>r_{\rm c}$, then the Cox point process associated to $\Lambda$ satisfies $\lcnp=0$, analogously to the case of Example~\ref{ex-mixedPPP}. In contrast, for $r<r_{\rm c}$ we have $\lcnp=\infty$. Note that it is also easy to see that that $\lcnp(r_{\rm c}) \in \{ 0, \infty \}$.
\end{example}

It is an interesting open question whether there exists an ergodic directing measure where for the associated Cox--Boolean model $\lcvp$, $\lcdp$ or $\lcnp$ are equal to zero without complete coverage, and also whether the non-existence of moments described in Part~\eqref{T3_NonSubMom} of Theorem~\ref{Thm_3} can occur for an ergodic directing measure $\Lambda$ without assuming $\E[\varrho^{d+s}]=\infty$.

\subsection{Deviations of the directing measure}\label{sec-overshoot}
In this section, we present examples for which the deviation condition, Condition~\eqref{Erg_cond}, holds. 
Let us start with the following statement for $b$-dependent environments with sufficient integrability.  
\begin{lemma}\label{lem_5}
Let $\L$ be stationary and $b$-dependent for some $b>0$. Then, Condition~\eqref{Erg_cond} holds if 
\begin{align*}
\E\big[\big|\L(B_b)-|B_b|\big|^{2k}\big]&<\infty\qquad\text{for some }(1+s)/d<k\in\N.
\end{align*}
\end{lemma}
The proof rests on an application of Lemma~\ref{lem_4} and is presented in the Section~\ref{sec-Exampleproofs}.  
Note that Lemma~\ref{lem_5} can be used to show that Condition~\eqref{Erg_cond} holds for example for the shot-noise field as presented in Example~\ref{Ex_2}. 
Indeed, in case $b={\rm diam}\big({\rm supp}(\kappa)\big)<\infty$ and $\kappa_{\rm max}=\sup\{\kappa(x)\colon x\in \R^d\}<\infty$, we have for any $n\in\N$
\begin{align*}
\E\big[\L(B_b)^n\big]\le (c_1|B_b|)^n\E\big[\#(Y\cap B_{2b})^n\big]\le c(\kappa_{\rm max}|B_b|)^n(\mu|B_{2b}|)^n, 
\end{align*}
where $c$ is a constant depending only on $n$ and $\mu$, the intensity of the underlying Poisson point process. But this implies that for all $k\in \N$
\begin{align*}
\E\big[|\L(B_b)-|B_b||^{2k}\big]<\infty. 
\end{align*}

Let us next verify that the unbounded stabilizing absolutely continuous directing measure exhibited in Example~\ref{Ex_3} also satisfies Condition~\eqref{Erg_cond}. For this, we verify Condition~\eqref{eq:cond2} in Lemma~\ref{lem_4} for $\beta=2k$ with $k\in\N$ sufficiently large. For the directing measure presented in Example~\ref{Ex_3} we can re-express the $2k$-th central moment in terms of the complete Bell polynomials $G_n$, see~\cite{GST19}, i.e., 
\begin{align*}
\E\big[\big(\L(B_\a)-|B_\a|\big)^{2k}\big]=\E\Big[\Big(\sum_{i\in I}|B_{\varrho_i}(Y_i)\cap B_\a|-|B_\a|\Big)^{2k}\Big]=G_{2k}(0,u_2,\dots,u_{2k}),
\end{align*}
where  
\begin{align*}
u_j=\int\int \big( \la\mu|B_{r}(x)\cap B_\a|\big)^{j}\nu(\d r)\d x\le (\la\mu)^j(2\a)^d\int r^{d(j+1)}\nu(\d r)
\end{align*}
is the $j$-th cumulant of $\sum_{i\in I}|B_{\varrho_i}(Y_i)\cap B_\a|$ and $\nu=\P\circ \varrho^{-1}$. In particular, using the recursive formula for Bell polynomials, if $\int r^{d(2k+1)}\nu(\d r)<\infty$, then there exists a constant $c>0$, only depending on $d,k,\nu,\la$ and $\mu$, such that 
\begin{align*}
G_{2k}(0,u_2,\dots,u_{2k})\le c\a^{dk}. 
\end{align*}
Hence, Condition~\eqref{eq:cond2} is verified once we have that $s<dk$, which is the case if $\int r^{d(2k+1)}\nu(\d r)<\infty$ for some $k>s/d$.  

Finally, we establish Condition~\eqref{Erg_cond} for the stationary singular but not $b$-dependent environment defined via the Poisson--Delaunay tessellation, as presented in Example~\ref{Ex_4}. 
\begin{lemma}\label{lem_6}
For any $s>0$, Condition~\eqref{Erg_cond} holds for $\L(\d x)=\nu_1(S\cap \d x)$, where $S$ is the edge set of the Delaunay tessellation based on a homogeneous Poisson point process $Y$ in $\R^2$. 
\end{lemma}
We give the proof in Section~\ref{sec-Exampleproofs}. Let us note that Lemma~\ref{lem_6} also holds for example for the Poisson--Voronoi tessellation in $\R^2$, but we do not present a proof for this statement.

\section{Proofs}\label{Sec-Proofs}
Recall for all $x\in\R^d$ the notations $Q_r(x)=Q_r+x$, where $Q_r$ is the box of side length $2r$ centered at the origin, and $B_r(x)=B_r+x$, where $B_r$ is the ball of radius $r$ centered at the origin. Further, let $\nu$ denote the distribution of $\varrho$ and $v_d=|B_1|$. Also, for any set $A\subset\R^d$, we define the Boolean model based on the points in $A$,
$$\C(A)=\bigcup_{i \in I\colon X_i\in A} B_{\varrho_i}(X_i),$$
and let $C_x(A)\subset\C(A)$ denote the cluster in $\C(A)$ that contains $x\in \R^d$. 
\subsection{Proofs of Propositions~\ref{Prop_1}, \ref{Prop_2}, and Theorem~\ref{Prop_3}}\label{sec-Basicproofs}
\begin{proof}[Proof of Proposition~\ref{Prop_1}]
Let $V_n=n^d- |\C\cap Q_{n/2}|$ denote the volume of the vacant region in $Q_{n/2}$. Then, we have the following estimate for the expected volume fraction,
\begin{align*}
n^{-d}\E[V_n]&=n^{-d}\E\Big[\int_{Q_{n/2}}\one\{x\not\in \C\}\d x\Big]=\P(o\not\in \C)\\
&=\E\Big[\prod_{i\in I}\P(|X_i|\ge \varrho_i)\Big]=\E\Big[\exp\Big(\sum_{i\in I}\log\P(|X_i|\ge \varrho_i)\Big)\Big]\\
&=\E\Big[\exp\Big(-\lambda\int\Lambda(\d x)\P(\varrho> |x|)\Big)\Big]=\E\Big[\exp\Big(-\lambda\E[\Lambda(B_\varrho)]\Big)\Big]\\
&\ge \exp\big(-\lambda\E[\Lambda(B_\varrho)]\big)=\exp\big(-\lambda v_d\E[\varrho^d]\big)>0,
\end{align*}
where we used first stationarity, then the representation of the Laplace transform for i.i.d.~marked Poisson point processes and finally Jensen's inequality. Now, $\E[V_n]>0$ implies that $\P(V_n>0)>0$ for all $n\in\N$ and hence with positive probability $\C\neq \R^d$. Since this is a translation-invariant event, under ergodicity, $\P(\C=\R^d)=0$. 
\end{proof}

\begin{proof}[Proof of Proposition~\ref{Prop_2}]
As in the proof of Proposition~\ref{Prop_1}, by Jensen's inequality,
$$\P(o\not\in \C)=\E[\exp(-\lambda\E[\Lambda(B_\varrho)])]\ge \exp(-\lambda v_d\E[\varrho^d])>0.$$
Since $\Lambda$ is ergodic, we conclude using the same arguments as in~\cite[Corollary 3.1]{MR96}. 
\end{proof}

\begin{proof}[Proof of Theorem~\ref{Prop_3}]
The proof rests on a generalization of the Burton--Keane argument presented in the proof of~\cite[Theorem 3.6]{MR96} additionally using the FKG inequality for Poisson point processes. 
First note that, by ergodicity, the number of infinite clusters in $\C$ is almost surely equal to a constant $K \in \N_0 \cup \{ \infty \}$ for all possible choices of the parameters, see for example~\cite[Section 3.6]{MR96}. Let us first exclude the case that $2\le K<\infty$. We do this by contradiction and assume that the number of unbounded clusters is equal to $K$ almost surely for some $K<\infty$. 
Let $r=\esssup(\varrho)<\infty$ and $\eps>0$ and consider the event
\begin{align*}
E_{n}=\{&\text{for all unbounded clusters }C\subset \C(Q_{3n}^{\rm c})\colon \dist(C, Q_{3n})<r-\eps\text{ and}\\
&\text{there exists }X_1,\dots,X_l\in X\cap Q_{3n}\text{ such that }|X_j-X_{j+1}|< \varrho_j+\varrho_{j+1}\text{ for all }j=1,\dots, l-1,\\
&B_{\varrho_1}(X_1)\cap C\neq\emptyset \text{ and }X_l\in Q_n\}. 
\end{align*}
In words, in the event $E_n$, all the unbounded clusters, restricted to the parts that come from Cox points outside $Q_{3n}$, intersect a small neighborhood around $Q_{3n}$ and are moreover connected to $Q_n$, with one point in $Q_n$. Note however that in the event $E_n$, the Cox points in $Q_{3n}$ that establish the connection to $Q_n$ do not necessarily have to be different for different unbounded clusters of $\C(Q_{3n}^{\rm c})$. Next, consider the event 
\begin{align*}
F_{n}=\{&\text{for all $x,y\in \supp(\L)\cap Q_n$ there exist }X_1,\dots,X_l\in X\cap Q_{2n}\text{ such that}\\
&|X_j-X_{j+1}|< \varrho_j+\varrho_{j+1}\text{ for all }j=1,\dots, l-1, \, x\in B_{\varrho_1}(X_1)\text{ and }y\in B_{\varrho_l}(X_l)\}
\end{align*}
that any pair of locations in $\supp(\L)\cap Q_n$ is connected by a finite path of Cox points in $Q_{2n}$ with overlapping balls. Now, in the event $E_n\cap F_n$, the $K$ unbounded clusters of $\C(Q_{3n}^{\rm c})$ are in fact connected in the Boolean model $\C$ and hence $K\in\{0,1\}$. Thus, it suffices to prove that $\P(E_n\cap F_n)>0$. In particular, it suffices to prove that $\P(E_n\cap F_n\cap G_n)>0$, where $G_{n}=\{\sup\{R_y\colon y\in Q_{2n} \cap\,  \Q^d\} < n/2\}$ is the event which guarantees that $\supp(\L_{Q_n})$  is $r$-connected in $\supp(\L_{Q_{2n}})$. For this, note that 
\begin{align*}
\P(E_n\cap F_n\cap G_n)&=\E[\one\{\L\in G_n\}\P(E_n\cap F_n |\L)]\ge \E[\one\{\L\in G_n\}\P(E_n|\L)\P(F_n |\L)],
\end{align*}
where we used the FKG inequality. To see that the FKG inequality is indeed applicable, for any measurable set $B \subseteq \R^d$, we write $P^\L_B$ for a Poisson point process with intensity measure $\la\L_B$. Then, 
\begin{align*}
\P(E_n\cap F_n |\L)=\int \int \one\{\omega_1\cup\omega_2\in E_n\}\one\{\omega_2\in F_n\}P^\L_{Q_{2n}}(\d\omega_2)P^\L_{Q_{2n}^{\rm c}}(\d\omega_1), 
\end{align*}
where we used the independence of the Poisson point processes and the fact that $F_n$ only depends on Poisson points in $Q_{2n}$. Conditioned on $\omega_1$, the functions $\omega_2\mapsto \one\{\omega_1\cup\omega_2\in E_n\}$ and $\omega_2\mapsto \one\{\omega_2\in F_n\}$ are both increasing, and hence, by the FKG inequality for Poisson point processes, see for example~\cite[Theorem 20.4]{LP17},
\begin{align*}
\P(E_n\cap F_n |\L)&\ge \int \int \one\{\omega_1\cup\omega_2\in E_n\}P^\L_{Q_{2n}}(\d\omega_2)P^\L_{Q_{2n}^{\rm c}}(\d\omega_1)\int \one\{\omega_2\in F_n\}P^\L_{Q_{2n}}(\d\omega_2)\\
&= \P(E_n|\L)\P(F_n |\L). 
\end{align*}
To continue, let us abbreviate $e_n(\L)=\P(E_n|\L)$ and $f_n(\L)=\P(F_n |\L)$. Note that the event $E_n$ is implied by the event that all $K$ unbounded clusters have a non-empty intersection with the box $Q_n$. Hence, by ergodicity, for all $\eps_1$, there exists $n_1\in\N$ such that for all $n\ge n_1$, $\P(E_n)>1-\eps_1$. Hence, there exists $c>0$ such that $\P\big(H_n\big)>1-\eps_1$, where $H_n=\{e_n(\cdot)>c\}$. This can be easily checked via contradiction. 
Moreover, by the assumption of essential $r$-connectedness, for all $\eps_2>0$, there exists $n_2\ge n_1$ such that for all $n\ge n_2$ we have that $\P(G_n)>1-\eps_2$ and hence, for sufficiently small $\eps_1$ and $\eps_2$, also $\P(G_n\cap H_n)>0$ for all sufficiently large $n$.
Hence, for such $n$,
\begin{align*}
\E[\one\{\L\in G_n\}e_n(\L)f_n(\L)]\ge \E[\one\{\L\in G_n\cap H_n\}e_n(\L)f_n(\L)]>c\E[\one\{\L\in G_n\cap H_n\}f_n(\L)]
\end{align*}
and, in order for the expression on the right-hand side to be positive, it suffices to prove that $f_n(\L)>0$ for almost all $\L\in G_n\cap H_n$. 
But this is true since 
conditioned on $\L$, under $G_n$, points can be connected via sufficiently many Poisson points, which has positive probability. 

\medskip
What remains to be proved is that $\P(K=\infty)=0$. Consider the event 
\begin{align*}
E'_n=\{&\text{there is an unbounded cluster $C$ that contains a Cox point } X_1\in Q_{n}  \text{ and }\\
&C\cap Q_{2n}^{\rm c} \text{ contains at least three unbounded clusters}\}.
\end{align*}
Using the exact same arguments as in the proof of~\cite[Theorem 3.6]{MR96}, the claim follows once we have shown that $\P(E'_n)>0$. In order to show this, recall the definitions of the events $F_n$ and $G_n$ from above 
and note that $\P(E'_n)\ge \P(E''_n\cap F_n\cap G_n)$, where 
\begin{align*}
E''_{n}=\{&\text{there exist at least three unbounded clusters }C\subset \C(Q_{3n}^{\rm c})\colon \dist(C, Q_{3n})<r-\eps\text{ and}\\
&\text{there exists }X_1,\dots,X_l\in X\cap Q_{3n}\text{ such that }|X_j-X_{j+1}|< \varrho_j+\varrho_{j+1}\text{ for all }j=1,\dots, l-1,\\
&B_{\varrho_1}(X_1)\cap C\neq\emptyset \text{ and }X_l\in Q_n\}. 
\end{align*}
But this probability is positive, which can be seen using the same arguments as in the first part of the proof. 
\end{proof}

\subsection{Proof of Theorem~\ref{Thm_3}}\label{sec-Mainproofs}
We prove the four parts of Theorem~\ref{Thm_3} individually. 
Let us start with the proof of Part~\eqref{T3_SupCrit} which uses coupling arguments.
\subsubsection{Existence of a supercritical phase for percolation}
 
\begin{proof}[Proof of Part~\eqref{T3_SupCrit} of Theorem~\ref{Thm_3}]
As for the assertion that $\lambda_{\rm vp}<\infty$ in case $\Lambda$ is $\phi$-stabilizing and $\esssup(\varrho)$ is sufficiently large, let us note that it was observed in \cite[Corollary 2.5]{T18} that if $\Lambda$ is $\phi$-stabilizing, then
\[ R=\inf \Big\{ r > 0 \colon \lcvp<\infty \text{ for the Boolean model $\bigcup_{i \in I} B_r(X_i) $} \Big\}  \]
is positive and finite. 
Hence, the proof of the assertion reduces to a simple coupling argument as follows. Let us assume that $\esssup(\varrho)>R$, then there exists $r>R$ such that $p=\P(\varrho>r)>0$. Now, for $i \in I$ let
\[ \varrho_i^{-} = r \mathds 1 \{ \varrho_i \geq r \} \]
and consider $X^-=\{ X_i \colon i \in I, \varrho_i^{-}=r \}$, which is an independent thinning of the Cox point process $X$ with probability $p$, and hence a Cox point process with intensity $p\la\L$ thanks to the Colouring Theorem, see~\cite{K93}. Further, since $\varrho_i^- \leq \varrho_i$ almost surely, we have $\C^-=\bigcup_{X_i \in X^{-}} B_{\varrho_i^-}(X_i) \subset \mathcal C$. Now, by the definition of $R$, one can choose $\lambda>0$ such that $\C^-$ exhibits volume percolation with probability 1, and hence so does $\mathcal C$. This together with the stationarity of $\mathcal C$ implies that $|C_o|=\infty$ holds with positive probability, and therefore $\lcvp<\infty$. Again, $\lcdp<\infty$, since $\lcdp\le\lcvp$. Finally, since in the supercritical regime for volume percolation, the unbounded cluster in $\C^-$ already contains infinitely many points, so does $\C$, which implies that also $\lcnp<\infty$.
\end{proof}

\subsubsection{Existence of a subcritical phase for percolation}
The proof of Part~\eqref{T3_SubCrit} of Theorem~\ref{Thm_3} rests on a generalization of arguments first presented in~\cite{G08}, for Poisson--Boolean models, which leverage scaling properties of the process. Let us mention again that, for the existence of subcritical regimes for volume percolation and diameter percolation, alternatively the conditions of~\cite[Theorem 2.8]{G09} can be verified. However, we decided to present here a self-contained direct proof since the verification of the conditions in~\cite{G09} is not substantially quicker, less instructive, and less transferable towards point-number percolation. 

Recall that $C_x(A)$ denotes the connected component of the Boolean model based on points in $A$ that contains $x$, then we define for any $x\in \R^d$ and $\a>0$ the event
\begin{align*}
G(x,\a)&=\{\big(C_x(B_{10\a}(x))\cup B_\a(x)\big)\not\subset B_{8\a}(x)\}
\end{align*}
that the cluster of $x$, only using points in $B_{10\a}(x)$, reaches beyond $B_{8\a}(x)$. Then, we have the following lemma.  
\begin{lemma}\label{Lem_3}
Consider the Boolean model based on the Cox point process $X$ with stationary intensity measure $\la\L$, where $\L$ is $\phi$-stabilizing. Then, there  exists a constant $c>0$, depending only on $d$, such that for all $\a>0$, we have 
\begin{align*}
\P(G(o,10\a))&\le c\P(G(o,\a))^2+\lambda c\int_\a^\infty r^d\nu(\d r)+2c\phi(10\a)\text{ and}\\
\P(G(o,\a))&\le c\lambda\a^d.
\end{align*}
\end{lemma}
\begin{proof}[Proof of Lemma~\ref{Lem_3}]
Consider the event
\begin{align*}
H'(\a)&=\{\text{there exists }X_i\in X\cap B_{100\a}\colon \varrho_i\ge \a\},
\end{align*}
and let  $K_\a$ denote a finite subset of the sphere $S_\a=\{x\in\R^d\colon |x|=\a\}$ such that 
$$S_\a\subset \bigcup_{x\in K_\a}B_1(x).$$
Then, a small generalization of the arguments presented in the proof of~\cite[Lemma 3.3]{G08} gives 
\begin{align*}
\P(G(o,10\a))&\le\P\Big(\Big( \bigcup_{k\in K_{10}} G(\a k,\a)\Big) \cap \Big(\bigcup_{l\in K_{100}}G(\a l,\a)\Big)\Big)+\P(H'(\a))\\
&\le\sum_{k\in K_{10},\, l\in K_{100}}\E[\P(G(\a k,\a)|\L)\P(G(\a l,\a)|\L)]+\P(H'(\a)),
\end{align*}
since $\big(G(o,10\a)\setminus H(\a)\big)\subset \Big(\Big(\bigcup_{k\in K_{10}} G(\a k,\a)\Big)\cap \Big(\bigcup_{l\in K_{100}}G(\a l,\a)\Big)\Big)$. Then, using $\phi$-stabilization, 
\begin{align*}
&\E\big[\P(G(\a k,\a)|\L)\P(G(\a l,\a)|\L)\big]\\
&\le \E\big[\one\{\sup_{y \in Q_{10\a}(\a k) \cap \, \Q^d}R_y < \a\}\P(G(\a k,\a)|\L)\one\{\sup_{y \in Q_{10\a}(\a l) \cap \, \Q^d}R_y < \a\}\P(G(\a l,\a)|\L)\big]+2\phi(10\a),
\end{align*}
and hence, for some $c_1>0$ only depending on the dimension, by stationarity,
\begin{align*}
\P(G(o,10\a))&\le  c_1\P(G(o,\a))^2+2c_1\phi(10\a)+\P(H'(\a)). 
\end{align*}

\medskip
Finally, as in~\cite[Lemmas 3.5 and 3.6]{G08}, there exist constants $c_2,c_3>0$, depending only on the dimension, such that 
\begin{align*}
\P(H'(\a))&\le \lambda c_2\int_\a^\infty r^d\nu(\d r)\qquad \text{ and}\qquad\P(G(o,\a))\le \lambda c_3\a^d.
\end{align*}
This finishes the proof. 
\end{proof}
Next, recall that by $X(A)$ we denote number of points of $X$ in $A\subset \R^d$. Consider
$$M=\sup\{|x|\colon x\in C_o\},$$
the largest distance to the origin of any point in the cluster of the origin. We have the following estimates. 
\begin{lemma}\label{Lem_4}
Consider the Boolean model based on the Cox point process $X$ with stationary intensity measure $\la\L$. Then, there  exists a constant $c>0$, only depending on $d$, such that for all $\a>0$, we have 
\begin{align*}
\P(M\ge 9\a)&\le \P(G(o,\a))+\lambda c\int_\a^\infty r^d\nu(\d r)\text{ and}\\ 
\P(X^\la(C_o)> X^\la(B_{8\a}))&\le \P(G(o,\a))+\lambda c\int_\a^\infty r^d\nu(\d r).
\end{align*}
\end{lemma}
\begin{proof}[Proof of Lemma~\ref{Lem_4}]
Consider the event, 
\begin{align*}
H(\a)&=\{\text{there exists }X_i\in X\cap B^{\rm c}_{10\a}\colon |X_i|\le 9\a+\varrho_i\}.
\end{align*}
Then, by~\cite[Lemma 3.2]{G08}, 
\begin{align*}
\P(M\ge 9\a)&\le \P(G(o,\a))+\P(H(\a)). 
\end{align*}
Moreover, we have that $\big(G(o,\a)^{\rm c}\cap H(\a)^{\rm c}\big)\subset \{C_o\subset B_{8\a}\}$ since, under the event $H(\a)^{\rm c}$, points outside $B_{9\a}$ cannot help the cluster $C_o$ to reach outside of $B_{8\a}$. But since $\{C_o\subset B_{8\a}\}\subset \{X(C_o)\le X(B_{8\a})\}$, we have 
\begin{align*}
\P(X(C_o)> X(B_{8\a}))&\le \P(G(o,\a))+\P(H(\a)). 
\end{align*}
Finally, by~\cite[Lemma 3.4]{G08}, there exist constants $c>0$, only depending on the dimension, such that 
\begin{align*}
\P(H(\a))&\le \lambda c\int_\a^\infty r^d\nu(\d r),
\end{align*}
which finishes the proof.
\end{proof}
We will also need the following essential result about convergence and integrability properties of functions satisfying some scaling inequality. 
\begin{lemma}[{\cite[Lemma 3.7]{G08}}]\label{Lem_5}
Let $f$ and $g$ be two bounded measurable functions from $[1,\infty]$ to $[0,\infty)$. Additionally, let $f$ be bounded by $1/2$ on $[1,10]$, $g$ be bounded by $1/4$ on $[1,\infty]$ and assume 
\begin{align*} 
f (\a) \le f (\a/10)^2 + g(\a),\qquad\text{for all }\a\ge 10.
\end{align*}
Then, $\lim_{\a\uparrow\infty}g(\a)=0$ implies that $\lim_{\a\uparrow\infty}f(\a)=0$. Moreover, if $s\in[0,\infty)$ is such that $\int_1^{\infty}\a^sg(\a)\d\a<\infty$, then $\int_1^{\infty}\a^sf(\a)\d\a<\infty$.
\end{lemma}
The first statement of Lemma~\ref{Lem_5} can be used to prove existence of a subcritical phase for percolation. 
\begin{proof}[Proof of part~\eqref{T3_SubCrit} of Theorem~\ref{Thm_3}]
We assume that $\E[\varrho^d]<\infty$ and $\L$ to be $\phi$-stabilizing. In order to prove that $\lcvp>0$ and $\lcdp>0$, it suffices to show that $\lim_{\a\uparrow\infty}\P(M\ge 9\a)=0$ for all sufficiently small $\la$. For $\lcnp>0$, it suffices to show that $\lim_{\a\uparrow\infty}\P(X^\la(C_o)> X^\la(B_{8\a}))=0$ for all sufficiently small $\la$. But those two statements are true if $\lim_{\a\uparrow\infty}\P(G(0,\a))=0$, by an application of Lemma~\ref{Lem_4}. In order to show $\lim_{\a\uparrow\infty}\P(G(0,\a))=0$ we apply the Lemmas~\ref{Lem_3} and~\ref{Lem_5} for proper choices of $f$ and $g$. For this, first define 
$$\a_{\rm c}=\inf\{x\ge 1\colon \phi(\a)<(4c)^{-2}\text{ for all }\a\ge x\},$$
and note that $\a_{\rm c}<\infty$ by the assumption of $\phi$-stabilization. Based on this, we make the following definitions. In case $\E[\varrho^d]\ge 1$,
\begin{align*}
A=10\a_{\rm c}\quad&\text{ and }\quad\lambda_o=\big(2 c^2 \E[\varrho^d](100\a_{\rm c})^d\big)^{-1}, \numberthis\label{Alambda0def} \\
f(\a)=c\P(G(o,A\a))\quad&\text{ and }\quad g(\a)=\lambda c^2\int_{A\a/10}^\infty r^d\nu(\d r)+2c^2\phi(A\a/10). \numberthis\label{fgdef}
\end{align*}
Then, indeed, using Lemma~\ref{Lem_3}, we have that 
\begin{align*}
f(\a)&\le c^2\lambda(A\a)^d\le \frac{(10\a_{\rm c}\a)^d}{2\E[\varrho^d](100\a_{\rm c})^d}\le \frac{1}{2},\qquad\text{ for all }1\le \a\le10\text{ and }\la<\la_o,\\
g(\a)&\le \frac{1}{8\a_{\rm c}^d}+2c^2\phi(\a_{\rm c}\a)\le \frac 1 4,\qquad\text{ for all }1\le \a\text{ and }\la<\la_o,\text{ and}\\
f(\a)&\le c^2\P(G(o,A\a/10))^2+ \lambda c^2\int_{A\a/10}^\infty r^d\nu(\d r)+2c^2\phi(A\a/10)=f(\a/10)^2+g(\a)\quad\text{ for all }\a \ge 10.
\end{align*}
Hence, since $\lim_{\a\uparrow\infty}g(\a)=0$, an application of Lemma~\ref{Lem_5} gives the result. 
On the other hand, in case $\E[\varrho^d]< 1$, we set
\begin{align*}
A=10\a_{\rm c}\qquad&\text{ and }\qquad\lambda_o=(2 c^2 (100\a_{\rm c})^d)^{-1},
\end{align*}
and we again define the functions $f$ and $g$ according to~\eqref{fgdef}. 
Then, again using Lemma~\ref{Lem_3}, we have that,
\begin{align*}
f(\a)&\le c^2\lambda(A\a)^d\le \frac{(10\a_{\rm c}\a)^d}{2(100\a_{\rm c})^d}\le \frac{1}{2},\qquad\text{ for all }1\le \a\le10\text{ and }\la<\la_o\\
g(\a)&\le \frac{1}{2(100\a_{\rm c})^d}+2c^2\phi(\a_{\rm c}\a)\le \frac 1 4,\qquad\text{ for all }1\le \a\text{ and }\la<\la_o,\text{ and}\\
f(\a)&\le c^2\P(G(o,A\a/10))^2+ \lambda c^2\int_{A\a/10}^\infty r^d\nu(\d r)+2c^2\phi(A\a/10)=f(\a/10)^2+g(\a)\quad\text{ for all }\a \ge 10.
\end{align*}
Hence, again an application of Lemma~\ref{Lem_5} gives the result. 
\end{proof}

\subsubsection{Existence of moments}
In this section we prove Part~\eqref{T3_SubMom} of Theorem~\ref{Thm_3} by establishing regimes of sufficiently small $\la$, such that moments for the volume, the diameter, and the number of points in the cluster of the origin exist. Again, let us mention that, for the existence of moments for cluster volumes and cluster diameters, alternatively the conditions of~\cite[Theorem 2.9]{G09} can be verified. However, we present here self-contained direct proofs based on the results in the previous section.  

\begin{proof}[Proof of Part~\eqref{T3_SubMom} of Theorem~\ref{Thm_3}]
In this proof we use the statement about integrability from Lemma~\ref{Lem_5}. Let $s>0$ and recall the definitions of the functions $f$ and $g$ from \eqref{fgdef}. 
Under the assumption that $\int_{0}^\infty\a^{s-1}\phi(\a)\d \a<\infty$ and $\E[\varrho^{d+s}] < \infty$ we have that  
\begin{align*}
\int_{0}^\infty \a^{s-1}g(\a)\d \a&=\lambda c^2\int_{0}^\infty \a^{s-1}\int_{\a_{\rm c}\a}^\infty r^d\nu(\d r)\d\a+2c^2\int_0^\infty\a^{s-1}\phi(\a_{\rm c}\a)\d \a\\
&=\lambda c^2\a_{\rm c}^{-s}\int_{0}^\infty \a^{s-1}\int_{\a}^\infty r^d\nu(\d r)\d\a+2c^2\a_{\rm c}^{-s}\int_0^\infty\a^{s-1}\phi(\a)\d \a\\
&=s^{-1} \lambda c^2\a_{\rm c}^{-s}\int_{0}^\infty r^{d+s}\nu(\d r)+2c^2\a_{\rm c}^{-s}\int_0^\infty\a^{s-1}\phi(\a)\d \a<\infty.
\end{align*}
Hence, for $\la$ sufficiently small, i.e., $\la<\la_o$ where $\la_o$ was defined in~\eqref{Alambda0def}, by an application of the Lemmas~\ref{Lem_5} and~\ref{Lem_4}, we obtain
\begin{align*}
\infty&>\int_{0}^\infty \a^{s-1}f(\a)\d \a=c\int_{0}^\infty \a^{s-1}\P(G(o,A\a))\d\a\\
&\ge c\int_{0}^\infty \a^{s-1}\P(M\ge 9A\a)\d \a-\lambda c^2\int_{0}^\infty \a^{s-1}\int_{A\a}^\infty r^d\nu(\d r)\d \a\\
&= c(9A)^{-s}\int_{0}^\infty \a^{s-1}\P(M\ge \a)\d \a-\lambda c^2A^{-s}\int_{0}^\infty r^{d+s}\nu(\d r)\\
&= c(9A)^{-s}\E[M^s]-\lambda c^2A^{-s}\E[\varrho^{d+s}]. 
\end{align*}
Hence, $\E[M^s]<\infty$. Since $\diam(C_o)\le M$ and $|C_o|\le M^d$, also $\E[\diam(C_o)^s]<\infty$ and $\E[|C_o|^{s/d}]<\infty$. This provides the proof for the existence of moments for the volume and diameter of $C_o$. 

\medskip
For the existence of moments for the number of points in $C_o$, note that under the assumptions $\int_{0}^\infty\a^{s-1}\phi(\a)\d \a<\infty$ and $\E[\varrho^{d+s}] < \infty$, we have, that 
\begin{equation}\label{Eq10}
\begin{split}
\E[X(C_o)^{s/d}]&=\int_{0}^\infty \a^{s-1}\P(X(C_o)\ge \a^d)\d \a\\
&\le \int_{0}^\infty \a^{s-1}\P(X(B_{K\a})\ge \a^d)\d \a+\int_{0}^\infty \a^{s-1}\P(\diam(C_o)\ge K\a^d)\d \a,
\end{split}
\end{equation}
for all $K>0$, where the second summand in~\eqref{Eq10} is finite for all sufficiently small $\lambda$, for any choice of $K>0$. Hence, by our assumptions, it suffices to show that the first summand in~\eqref{Eq10} is finite for all sufficiently small $\lambda$, which is equivalent to showing that 
\begin{align*}
\int_{0}^\infty \a^{s-1}\P(X(B_{\a})\ge c'\a^d)\d \a<\infty,
\end{align*}
for some $c'>0$ and $\lambda>0$. But this is true whenever $c'>c$ for any $c>0$ for which the Assumption~\eqref{Erg_cond} is satisfied for the constant $c/\la$. Indeed, for these choices, 
\begin{align*}
\P(X(B_{\a})\ge c'\a^d)\le \E[\one\{\la\L(B_\a)< c\a^d\}\P(X(B_{\a})\ge c'\a^d|\L)]+\P(\la\L(B_\a)\ge c\a^d), 
\end{align*}
where by assumption $\int_{0}^\infty \a^{s-1}\P(\L(B_\a)\ge c\a^d/\la)\d \a<\infty$. 
For the other term, we can apply Poisson concentration inequalities, which lead to an upper bound given by
\begin{align*}
\E[\one\{\lambda\L(B_\a)< c\a^d\}\P(X(B_{\a})\ge c'\a^d|\L)]&\le\exp(-c''\a^d),
\end{align*}
where $c''=c\big((c'/c)\log(c'/c)+1-c'/c\big)>0$, and hence 
\begin{align*}
\int_{0}^\infty \a^{s-1}\exp(-c''\a^d)\d \a<\infty. 
\end{align*}
This concludes the proof. 
\end{proof}

\subsubsection{Non-existence of moments}
In this part we finish the proof of Theorem~\ref{Thm_3} by verifying Part~\eqref{T3_NonSubMom} using a generalization of the proof presented for~\cite[Theorem 3.2]{MR96}.

\begin{proof}[Proof of Part~\eqref{T3_NonSubMom} of Theorem~\ref{Thm_3}]
Let $\E[\varrho^{d+s}]=\infty$ and assume $\L$ to be ergodic. We may assume that $\E[\varrho^{d}]<\infty$, otherwise there is nothing to prove, see Theorem~\ref{Thm_1} Part~\eqref{co_co}. Let us start by considering moments of the diameter of $C_o$. We reproduce the strategy used for~\cite[Lemma 3.9]{G08}. Note that for all $\a\ge 0$
\begin{align*}
\P(\diam(C_o)\ge \a)&\ge \P(\text{there exists }X_i\in X\text{ such that }|X_i|+\a < \varrho_i)\\
&= \E\Big[1-\exp\big(-\la\int_0^\infty\L(B_{r+\a})\nu(\d r)\big)\Big]= 1- \E\Big[\exp\big(-\la\int_\a^\infty\L(B_r)\nu(\d r)\big)\Big].
\end{align*}
By ergodicity of $\L$, for all $\eps \in (0,1)$, there exists an $R>0$ such that the event 
\begin{align*}
\{(1-\eps)|B_{r}|<\L(B_{r}), \text{ for all }r>R\}
\end{align*}
has probability at least $1/2$, and hence, conditioned on this event, for $\a >R$, we can further bound
\begin{align*}
\P(\diam(C_o)\ge \a)&\ge \frac 1 2- \frac 1 2 \E\Big[\exp\big(-\la v_d(1-\eps)\int_\a^\infty r^d\nu(\d r)\big)\Big].
\end{align*}
Since $\E[\varrho^{d}]<\infty$, there exists a constant $c>0$ such that 
\begin{align*}
\frac 1 2- \frac 1 2 \E\Big[\exp\big(-\la v_d(1-\eps)\int_\a^\infty r^d\nu(\d r)\big)\Big]\ge c\int_\a^\infty r^d\nu(\d r),
\end{align*}
and thus, since $\E[\varrho^{d+s}]=\infty$,
\begin{align*}
\E[\diam(C_o)^{s}]&\ge \int_R^\infty \a^{s-1}\P(\diam(C_o)\ge \a)\d \a\ge c\int_R^\infty \a^{s-1}\int_\a^\infty r^d\nu(\d r)\d\a=\infty.  
\end{align*}
This proves the case for the diameter. For the volume, we can use the same arguments. In particular, for all $\a>R$, there exists a finite constant $c>0$ such that we can estimate 
\begin{align*}
\P(|C_o|\ge |B_\a|)&\ge \P(\text{there exists }X_i\in X\text{ such that }|X_i|+\a < \varrho_i)\ge c\int_\a^\infty r^d\nu(\d r),
\end{align*}
and thus, 
\begin{align*}
\E[|C_o|^{s/d}]&\ge \int_R^\infty \a^{s-1}\P(|C_o|\ge \a^d)\d \a= v_d^{s/d}\int_{R/v_d^{1/d}}^\infty \a^{s-1}\P(|C_o|\ge|B_\a|)\d \a
=\infty.  
\end{align*}
This proves the case for volumes. Finally, for the number of points in $C_o$, as in the proof for the existence of a subcrititical regime, we need extra arguments. We follow the general approach used in~\cite{MR96}.  It will be convenient to assume the radii to be integer-valued by setting $\rho=\lfloor \varrho\rfloor$. This is no restriction since $\E[\varrho^{d+s}]=\infty$ if and only if $\E[\rho^{d+s}]=\infty$. Note that conditioned on $\L$, the Poisson point process $X^\L$ can be seen as a superposition of independent Poisson point processes $X^{\L,j}$ with intensities $\P(\rho=j)\lambda\L(\d x)$. The number of direct neighbors $X_i$ of the origin in $C_o$ coming from the process $X^{\L,j}$ such that $|X_i|\le j/2$ is denoted by $N^\L_{j}=X^{\L,j}(B_{j/2})$. 
Note that here, different from our guiding example for this part of the proof,~\cite[Theorem 3.2]{MR96}, it will be very useful to have the extra factor $1/2$ in the radius. Further, let 
\begin{align*}
M^\L=\max\{j\ge 0\colon N^\L_{j}>0\},
\end{align*}
and put $M^\L=-1$ in case $N^\L_{j}=0$ for all $j\ge 0$. Note that the event $\{M^\L=m\}$ depends only on the processes $X^{\L,j}$ for $j\ge m$, where we write $X^{\Lambda,-1}=\emptyset$. Hence, 
\begin{align*}
\E[X(C_o)^{s/d}|\L]&=\sum_{m=-1}^\infty\E[X(C_o)^{s/d}\one\{M^\L=m\}|\L]=\sum_{m=-1}^\infty\P(M^\L=m|\L)\E[X(C_o)^{s/d}|\L, M^\L=m]\\
&\ge \sum_{m>k_o}^\infty\P(M^\L=m|\L)\E[X(C_o)^{s/d}|\L, M^\L=m],
\end{align*}
where $k_o=\min\{i\ge 1\colon \P(\rho=i)>0\}$. Now, for the first term,
\begin{align*}
\P(M^\L=m|\L)&=\big(1-\P(N^\L_{m}=0|\L)\big)\prod_{j> m}\P(N^\L_{j}=0|\L)\\
&=\big(1-\exp\big(-\P(\rho=m)\lambda\L(B_{m/2})\big)\big)\exp\big(-\lambda\sum_{j> m}\P(\rho=j)\L(B_{j/2})\big)\\
&\ge\P(\rho=m) \lambda\L(B_{m/2})\exp\big(-\lambda\sum_{j\ge m}\P(\rho=j)\L(B_{j/2})\big), 
\end{align*}
where we used that $1-\exp(-x)\ge x\exp(-x)$ for all $x\ge 0$. 
For the second term, note that under the event $\{M^\L=m\}$ with $m> k_o$, there exists at least one Cox point $X_i$ with $|X_i|\le m/2$ and $B_m(X_i)\subset C_o$. In particular, $B_{m/2}\subset  B_m(X_i)$ and hence, for all Cox points $X_j\in B_{m/2}$ we have $X_j\in C_o$. Thus, we can estimate using independence, 
\begin{align*}
\E[X(C_o)^{s/d}|\L, M^\L=m]&\ge \E[X^{\L,k_o}(B_{m/2})^{s/d}|\L]. 
\end{align*}
Next, we distinguish two cases. 

{\bf Case $s/d\ge 1$:} In this case, $x\mapsto x^{s/d}$ is convex and hence, using Jensen's inequality,  
\begin{align*}
\E[X^{\L,k_o}(B_{m/2})^{s/d}|\L]\ge \big(\P(\rho=k_o)\lambda\L(B_{m/2})\big)^{s/d}. 
\end{align*}
Then, putting everything together, we have that 
\begin{align*}
\E[X(C_o)^{s/d}]&\ge \lambda^{1+s/d}\P(\rho_o=k_o)^{s/d} \sum_{m> k_o}\P(\rho=m)\E\big[\L(B_{m/2})^{1+s/d}\exp\big(-2\lambda\sum_{j\ge m}\L(B_{j/2})\P(\rho=j)\big)\big].
\end{align*}
To finish this case, by ergodicity of $\L$, there exists an $N\in\N$ such that the event 
\begin{align}\label{Event_1}
\{1-\eps<\frac{\L(B_{m/4})}{|B_{m/4}|}<1+\eps, \text{ for all }m>N\}
\end{align}
has probability at least $1/2$, and hence, conditioned on this event and assuming that $\E[\rho^d]<\infty$, 
\begin{align*}
\E[&X(C_o)^{s/d}]\ge \frac 1 2((1-\eps)\lambda)^{1+s/d}\\
&\times \P(\rho_o=k_o)^{s/d} \sum_{m\ge N\vee (k_o+1)}\P(\rho=m)|B_{m/4}|^{1+s/d}\exp\big(-2\lambda\sum_{j\ge m}|B_{j/2}|(1+\eps)\P(\rho=j)\big)=\infty,
\end{align*}
by the assumptions. Finally, we consider the other case. 
%

{\bf Case $s/d< 1$:} In this case, note that for a Poisson random variable $L$ with parameter $\mu$, we can estimate
\begin{align*}
\E[L^{s/d}]=\e^{-\mu}\sum_{n\ge 0}n^{s/d}\frac{\mu^n}{n!}=\mu\e^{-\mu}\sum_{n\ge 1}n^{s/d-1}\frac{\mu^{n-1}}{(n-1)!}=\mu\e^{-\mu}\sum_{n\ge 0}(n+1)^{s/d-1}\frac{\mu^{n}}{n!}, 
\end{align*}
where $x\mapsto (x+1)^{s/d-1}$ is convex. Hence, again via Jensen's inequality, 
\begin{align*}
\E[L^{s/d}]\ge \mu(\mu+1)^{s/d-1}\ge \frac{\mu}{\mu+1}\mu^{s/d}. 
\end{align*}
We can now use this to estimate
\begin{align*}
\E[X^{\L,k_o}(B_{m/2})^{s/d}|\L]\ge \frac{c\L(B_{m/2})}{c\L(B_{m/2})+1}\big(c\L(B_{m/2})\big)^{s/d},
\end{align*}
where we abbreviated $c=\P(\rho=k_o)\lambda$. Putting things together and using again ergodicity via the event~\eqref{Event_1}, we have
\begin{align*}
\E[X(C_o)^{s/d}]&\ge \lambda c^{s/d} \sum_{m> k_o}\P(\rho=m)\E\big[\L(B_{m/2})^{1+s/d}\frac{c\L(B_{m/2})}{c\L(B_{m/2})+1}\exp\big(-2\lambda\sum_{j\ge m}\L(B_{j/2})\P(\rho=j)\big)\big]\\
&\ge \frac 1 2\lambda c^{s/d}(1-\eps)^{1+s/d}\sum_{m\ge N\vee (k_o+1)}\P(\rho=m)|B_{m/4}|^{1+s/d}\frac{c|B_{m/2}|(1-\eps)}{c|B_{m/2}|(1+\eps)+1}\\
&\qquad\times\exp\big(-2\lambda\sum_{j\ge m}|B_{j/2}|(1+\eps)\P(\rho=j)\big)\big]. 
\end{align*}
Now, assuming that $\eps<1/3$, there exists $K\in\N$, such that $c|B_{m/2}|(1-\eps)/(c|B_{m/2}|(1+\eps)+1)> 1/2$, for all $m\ge K$. Hence, since $\E[\rho^d]<\infty$, we finally have that, 
\begin{align*}
\E[X(C_o)^{s/d}]&\ge \frac 1 4\lambda c^{s/d}(1-\eps)^{1+s/d}\sum_{m\ge N\vee K\vee (k_o+1)}\P(\rho=m)|B_{m/4}|^{1+s/d}\\
&\qquad\times \exp\big(-2\lambda\sum_{j\ge m}|B_{j/2}|(1+\eps)\P(\rho=j)\big)\big]=\infty. 
\end{align*}
This finishes the proof. 
\end{proof}

\subsection{Proofs of Lemma~\ref{lem_4}, Corollary~\ref{cor-bdd}, and Lemmas~\ref{lem_5} and~\ref{lem_6}}\label{sec-Exampleproofs}

\begin{proof}[Proof of Lemma~\ref{lem_4}]
The proof relies on two applications of the Markov inequality. For Condition~\eqref{eq:cond1}, using the exponential Markov inequality for some constant $\beta>0$, we have  
\begin{align}\label{Eq11}
\int_{0}^\infty \a^{s-1}\P(\L(B_\a)\ge c\a^d)\d \a\le \int_{0}^\infty \a^{s-1}\exp(-\beta c \a ^d)\E[\exp(\beta\L(B_\a))]\d \a.
\end{align}
By Condition~\eqref{eq:cond1}, there exists $\a_{\rm c},\beta>0$ and $c'>0$ such that 
\begin{align*}
\E[\exp(\beta\L(B_\a))]\le \exp(c'\a ^d),
\end{align*}
for all $\a>\a_{\rm c}$. Hence, for $c>c'/\beta$, the integral~\eqref{Eq11} is finite. 

For Condition~\eqref{eq:cond2}, writing $v_d=|B_1|$, choosing $c>v_d$ and using the Markov inequality with $x\mapsto x^\beta$, we have 
\begin{align*}
\int_{0}^\infty \a^{s-1}\P(\L(B_\a)\ge c\a^d)\d \a=(c-v_d)^\beta\int_{0}^\infty \a^{s-1}\a^{-\beta d}\E(|\L(B_\a)-|B_\a||^\beta)\d \a
\end{align*}
By Condition~\eqref{eq:cond2}, there exists $\a_{\rm c}>0$, $\beta>1$ and $c'>0$ such that 
\begin{align*}
\a^{s-1-\beta d}\E(|\L(B_\a)-|B_\a||^\beta)\le c'\a^{-1-\eps}
\end{align*}
for all $\a>\a_{\rm c}$, which shows the desired integrability.  
\end{proof}

\begin{proof}[Proof of Corollary~\ref{cor-bdd}]
There exists a coupling such that the Boolean model $\mathcal C$ described in the corollary is almost surely included in a Boolean model $\mathcal C'$ based on a Poisson point process with intensity $\lambda M \d x$ with the same radius distribution $\varrho$, where $M=\esssup( \ell_o)$ is the bounding constant. Hence, whenever the latter Boolean model satisfies $\lcvp>0$, $\lcdp>0$ or $\lcnp>0$, then the same assertion holds for $\mathcal C$. Similarly, for $s>0$, if $\mathcal C'$ satisfies $\lcve(s)>0$, $\lcde(s)>0$ or $\lcne(s)>0$, then the same holds for $\mathcal C$. We conclude that $\E[\varrho^d]<\infty$ implies all the assertions $\lcvp>0$, $\lcdp>0$ and $\lcnp>0$ for $\mathcal C$. Further, if for $s>0$, $\E[\varrho^{d+s}]<\infty$, then $\lcve(s/d)>0$, $\lcde(s)>0$ and $\lcne(s/d)>0$ hold for $\mathcal C$. 

Assume now that $\varrho$ is such that at least one of the assertions $\lcvp>0$, $\lcdp>0$ or $\lcnp>0$ fails for $\mathcal C'$. Then in fact all of these assertions fail and we have that $\E[\varrho^d]=\infty$. Hence, by Theorem~\ref{Thm_1} Part~\eqref{co_co}, we conclude that $\P(\mathcal C=\R^d)=1$, and hence $\lcvp=\lcdp=\lcnp=0$ for $\mathcal C$ as well, as required.

Finally, assume now that for some $s>0$ we have $\lcve(s/d)=0$ for $\mathcal C'$. Then it follows from Theorem~\ref{Thm_2} that $\E[\varrho^{d+s}]=\infty$. But then, since $\Lambda$ is ergodic, it follows from Theorem~\ref{Thm_3} that $\lcve(s/d)=0$ holds also for $\mathcal C'$. Using analogous arguments, we conclude that $\lcde(s)=\lcne(s/d)=0$ for $\mathcal C'$ implies $\E[\varrho^{d+s}]=\infty$ and thus also $\lcve(s)=\lcne(s/d)=0$ for $\mathcal C$. Therefore, the corollary follows.
\end{proof}

\begin{proof}[Proof of Lemma~\ref{lem_5}]
We verify Condition~\eqref{eq:cond2} from Lemma~\ref{lem_4} for $\beta=2k$. 
For convenience, let us replace $B_\a$ by cubes $Q_{\a/2}(x)$ of sidelength $\a>0$. Then, assuming $\a/b\in \N$, 
\begin{align*}
&\E\big[\big|\L(Q_{\a/2})-\a^d\big|^{2k}\big]=\E\big[\big|\sum_{z\in b\Z^d\cap Q_{\a/2} }Y^b_z\big|^{2k}\big]=\sum_{z_1,\dots, z_{2k}\in b\Z^d\cap Q_{\a/2} }\E\big[Y^b_{z_1}\cdots Y^b_{z_{2k}}\big],
\end{align*}
where $Y^b_z=\L(Q_{b/2}(z))-b^d$. Note that, as soon as one of the $z_i$ is isolated within $(z_1,\dots,z_{2k})$, we have $\E\big[Y_{z_1}\cdots Y_{z_{2k}}\big]=0$. But, the number of possible configurations of the $(z_1,\dots,z_{2k})$, such that none of the points is isolated is upper bounded by $(\a/b)^{dk}(2k(2d+1))^{k}$. Hence, using H\" older's inequality,
\begin{align*}
&\sum_{z_1,\dots, z_{2k}\in b\Z^d\cap Q_\a }\E\big[Y^b_{z_1}\cdots Y^b_{z_{2k}}\big]\le (\a/b)^{dk}(2k(2d+1))^{k}\E\big[(Y^b_o)^{2k}\big]. 
\end{align*}
Finally, under our assumptions, we have that $\E\big[(Y^b_o)^{2k}\big]<\infty$ and $k>1+(1+s)/d$, and thus $\limsup_{\a\uparrow\infty}\a^{s-d(2k-1)+\eps +dk}=0$ for some $\eps>0$. Hence, Condition~\eqref{eq:cond2} is satisfied. 
\end{proof}
\begin{proof}[Proof of Lemma~\ref{lem_6}]
First, we can use the Chebyshev inequality 
\begin{align*}
\int_{0}^\infty \a^{s-1}\P(\L(B_{\a})\ge c\a^2)\d\a\le c_1+ \int_{1}^\infty \a^{s-5}\E\big[\L(B_{\a})^2\big]\d\a.
\end{align*}
Next, we introduce the stabilization radii as in the proof of~\cite[Proposition 2.3]{JT19}. We define
\[ R=\min \{ r \in \N \colon r \geq 2\a\text{ and } \forall z \in \Z^d \text{ with } \Vert z \Vert_{\infty} = 2, Q_r(rz) \cap Y\neq \emptyset \}, \numberthis\label{NdefPDT}\]
the finest discretization of $\R^d$ into boxes such that every box in the $2$-annulus contains Poisson points in $Y$. Note that $R$ is almost surely finite. Then, for $k \in \N$ such that $k > \lceil 2\a \rceil$,
\[ \begin{aligned}
\P(R \geq k) & \leq \P\big( \exists z \in \Z^2 \text{ with } \Vert z \Vert_{\infty} =2 \text{ such that } Q_{k-1}((k-1)z) \cap Y  = \emptyset \big) 
 \leq 16 \exp(-\lambda(k-1)^2).
\end{aligned} \]
Note that once $k>\lceil 2\a \rceil$, the right-hand side does not depend on $\a$. 
%
Now, the crucial observation is that, on the event that $R=k$, Delaunay edges intersecting $B_\a$ must have both endpoints within the ball $B_{2k}$, for details see the proof of~\cite[Proposition 2.3]{JT19}. 
Hence, we can use H\" older's inequality
to estimate
\[ 
\begin{aligned}
\E\big[\L(B_{\a})^2\big]&\leq (2\a)^2\sum_{k \geq 2\a} \E \big[(Y(Q_{6k})^2 \mathds 1 \{ R = k \} \big] \\ 
& \leq  (2\a)^2\sum_{k \geq 2\a} \E \big[Y(Q_{6k})^4 \big]^{1/2} \P( R = k )^{1/2} \\
& \leq c_3 \a^2\sum_{k \geq 2\a } (36k^2\lambda)^2\exp(-\lambda(k-1)^2/2),
\end{aligned}
\]
for some constant $c_3>0$, coming from lower-order terms in the evaluation of the moments of the Poisson random variable $Y(B_{6k})$.  Hence, there exists a finite constant $c_4>0$ such that 
\begin{align*}
\int_{1}^\infty \a^{s-5}\E\big[\L(B_{\a})^2\big]\d\a\le c_4\sum_{k\ge1}k^{s+2}\exp(-\lambda(k-1)^2/2), 
\end{align*}
which is finite for any $s>0$. 
This finishes the proof. 
\end{proof}



\section{Acknowledgement}
The authors thank A.~Hinsen, C.~Hirsch and W.~K\"onig for interesting discussions and comments. 
The authors also thank an anonymous reviewer for suggesting (i) the reference~\cite{G09}, which provides an alternative proof strategy for the Part~\eqref{T3_SubCrit} and the first part of Part~\eqref{T3_SubMom} of Theorem~\ref{Thm_3}, and (ii) an alternative approach to the proof of the second part of Part~\eqref{T3_SubMom} of Theorem~\ref{Thm_3} via a simplified version of Condition~\eqref{Erg_cond}.
This work was funded by the Deutsche Forschungsgemeinschaft (DFG, German Research Foundation) under Germany's Excellence Strategy MATH+: The Berlin Mathematics Research Center, EXC-2046/1 project ID: 390685689 and by Orange S.A..


\end{document}